\documentclass[12pt]{article}

\usepackage{amsmath,amsbsy,amsfonts,amssymb,amsthm,dsfont,fullpage,units,bm}
\usepackage{graphicx}
\usepackage{subcaption}
\usepackage{algorithm,algorithmic,mathtools}
\usepackage{color,cases}
\usepackage{tikz}
\usepackage{multirow}
\usepackage{hyperref, todonotes}
\usepackage{authblk}
\usepackage{amsmath,amsbsy,amsfonts,amssymb,
amsthm,dsfont,fullpage,units,bm,mathrsfs,bbm,
empheq,graphicx,subcaption,algorithm,algorithmic,
mathtools,color,cases,tikz,multirow,hyperref,
comment,enumerate}

\hypersetup{linkcolor=black}

\hypersetup{
colorlinks = true,
citecolor = {blue},
citebordercolor = {white}
}

\usepackage[top=1.in, bottom=1in, left=1.1in, right=1.1in]{geometry}
\usepackage{bbding}
\usepackage{subcaption}
\usepackage{stmaryrd}
\usepackage{mathrsfs, booktabs, multirow}
\usepackage{amsmath}    
\usepackage{amscd}          
\usepackage{amssymb}
\usepackage{graphicx}
\usepackage{array}  
\usepackage{booktabs}       
\usepackage{xcolor}
\usepackage{colortbl}       
\usepackage{multicol}
\usepackage{amsthm}

\theoremstyle{plain}
\begingroup
\newtheorem{theorem}{Theorem}
\newtheorem*{theorem*}{Theorem}
\newtheorem*{"theorem"}{``Theorem''}

\endgroup

\theoremstyle{definition}
\begingroup

\endgroup

\theoremstyle{remark}
\begingroup

\endgroup

\DeclareMathOperator{\rad}{Rad}

\newcommand{\bi}{\begin{itemize}}
\newcommand{\ei}{\end{itemize}}

\providecommand{\R}{{\ensuremath{\mathbb{R}}}}

\newcommand{\beq}{\begin{equation}}
\newcommand{\eeq}{\end{equation}}
\newcommand{\beqn}{\begin{eqnarray}}
\newcommand{\eeqn}{\end{eqnarray}}
\newcommand{\bsub}{\begin{subequations}}
\newcommand{\esub}{\end{subequations}}
\newcommand{\bpm}{\begin{pmatrix}}
\newcommand{\epm}{\end{pmatrix}}

\newcommand{\CR}{\mathcal{R}}

\newcommand{\argmin}{\text{argmin }}
\newcommand\Def{\stackrel{\textrm{def}}{=}}

\newcommand{\cH}{\mathcal{H}}
\newcommand{\cL}{\mathcal{L}}

\newcommand{\cF}{\mathcal{F}}

\newcommand{\EE}{{\mathbb{E}}}
\newcommand{\E}{{\mathbb{E}}}

\newcommand{\RR}{\mathbb{R}}
\newcommand{\N}{\mathbb{N}}

\newcommand{\bom}{\bm{\omega}}
\newcommand{\ba}{\bm{a}}

\newcommand{\vb}{\bm{v}}

\newcommand{\bz}{{\bm{z}}}

\newcommand{\wb}{{\bm{w}}}
\newcommand{\xb}{{\bm{x}}}
\newcommand{\bu}{{\bm{u}}}

\newcommand{\bx}{{\bm{x}}}

\newcommand{\bv}{{\bm{v}}}
\newcommand{\bw}{\bm{w}}

\newcommand{\bF}{{\bm{F}}}

\newcommand{\bW}{{\bm{W}}}

\newcommand{\bU}{{\bm{U}}}
\newcommand{\bV}{{\bm{V}}}

\newcommand{\argmax}{\text{argmax }}

\newcommand{\bG}{\bm{G}}

\newcommand{\veps}{\varepsilon}

\newcommand{\dt}{\partial_t}
\newcommand{\dx}{\nabla_{\bm{x}}}

\newcommand*\diff{\mathop{}\!\mathrm{d}}

\begin{document}

\begin{center}
{\Large \bf Machine Learning and Computational Mathematics}

\vspace{.1in}

Weinan E~\footnote{\texttt{weinan@math.princeton.edu}}


\vspace{.1in}
{Princeton University and Beijing Institute of Big Data Research}

\end{center}
\centerline{\it In memory of Professor Feng Kang (1920-1993)}

\tableofcontents

\begin{abstract}
Neural network-based machine learning is capable of 
 approximating functions in very high dimension with unprecedented efficiency  and accuracy.
 This has opened up many exciting new possibilities, not just in traditional areas of artificial intelligence, but also in 
 scientific computing and computational science.  At the same time, machine learning has also acquired the reputation of 
  being a set of ``black box'' type of tricks,
 without fundamental principles. This has been a real obstacle for making further progress in machine learning.
 
In this article, we try to address the following two very important questions: 
(1) How machine learning has already impacted and will further impact computational mathematics, scientific computing and computational science?
(2) How computational mathematics, particularly numerical analysis, {can} impact machine learning?
We describe some of the most important progress that has been made on these issues.
Our hope is to put things into a perspective that will help to integrate machine learning with computational mathematics.

\end{abstract}

\section{Introduction}

Neural network-based machine learning (ML) has shown very impressive success on a variety of tasks in traditional
artificial intelligence. This includes classifying images, generating new  images such as (fake) human faces
and playing sophisticated games such as Go.  A common feature of all these tasks is that they involve objects in very 
high dimension. Indeed when formulated in mathematical terms,  the image classification problem is a problem of approximating
a high dimensional  function, defined on the space of images, to the discrete set of values corresponding to the
category of each image. The dimensionality of the input space is typically  3 times the number of pixels in the image, 
where 3 is the dimensionality of the color space. 
The image generation problem is a problem of generating samples from an unknown high dimensional distribution, given
a set of samples from that distribution.
The Go game problem is about solving a Bellman-like equation in dynamic programming, since the optimal strategy
satisfies such an equation. For sophisticated games such as Go, this Bellman-like equation is formulated on a huge space.

 All these are made possible by the ability to accurately approximate high dimensional functions,
 using modern machine learning techniques.
This opens up new possibilities for attacking problems that suffer from the ``curse of dimensionality'' (CoD):
As dimensionality grows, computational cost grows exponentially fast. This CoD  problem has been an essential obstacle for
the scientific community for a very long time.

Take, for example, the problem of  solving partial differential equations (PDEs) numerically.
With traditional numerical methods such as finite difference, finite element and spectral methods, we can now routinely
solve PDEs in three spatial dimensions plus the temporal dimension.  Most of the PDEs currently studied in computational  
mathematics belong to this category. Well known examples include
the  Poisson equation, the Maxwell equation, the Euler equation, the Navier-Stokes equations, and the PDEs for linear elasticity.
Sparse grids can increase our ability to handling PDEs to, say 8 to 10 dimensions.
This allows us to try solving problems such as the Boltzmann equation for simple molecules.
But we are totally lost when faced with PDEs, say in 100 dimension.
This makes it essentially impossible to solve  Fokker-Planck or Boltzmann equations for complex molecules, many-body Schr\"odinger, 
or the Hamilton-Jacobi-Bellman equations for realistic control problems.

This is exactly where machine learning can help. Indeed, starting with the work in \cite{ HanE2016, EHanJentzen2017, HanJentzenE2018},
machine learning-based numerical algorithm for solving high dimensional PDEs and control problems
has been one of the most exciting new developments
in recent years in scientific computing, and this has opened up a host of new possibilities for
computational mathematics.  We refer to \cite{Review-1} for a review of this exciting development.

Solving PDEs is just the tip of the iceberg.  There are many other problems for which CoD is the main obstacle, including:
\begin{itemize}
\item classical many-body problem, e.g. protein folding
\item turbulence. Even though turbulence can be modeled by the three dimensional Navier-Stokes equation, it has so many
active degrees of freedom that an effective model for turbulence should involve many variables. 
\item solid mechanics.  In solid mechanics, we do not even have the analog of the Navier-Stokes equation. Why is this the case?
Well, the real reason is that the behavior of solids is essentially a multi-scale problem that involves scales from atomistic all
the way to macroscopic.
\item multi-scale modeling. In fact most multi-scale problems for which there is no separation of scales belong to this category.
 An immediate example is the dynamics of polymer fluids or polymer melts.
\end{itemize}

Can machine learning help for these problems? More generally,
{\it can we extend the success of machine learning  beyond traditional AI?}
We will try to convince the reader that  this is indeed the case for many problems.

Besides being extremely powerful, neural network-based machine learning has also got the reputation of being a set of
tricks instead of a set of systematic scientific principles.  Its performance depends sensitively on the value of the hyper-parameters,
such as the network widths and depths, the initialization, the learning rates, etc. Indeed just a few years ago,
parameter tuning was considered to be very much of an art.  Even now, This is still the case for some tasks.
Therefore a natural question is:  Can we  understand these subtleties and  propose
better machine learning models whose performance is more robust?

In this article, we review what has been learned on these two issues.  We discuss the impact that machine learning has already made
or will make on computational mathematics, and how the ideas from computational mathematics, particularly numerical analysis,
can be used to help understanding and better formulating machine learning models.
On the former, we will mainly discuss the new problems that can now be addressed using ML-based algorithms.
Even though machine learning also suggests new ways to solve some traditional problems in computational mathematics,
we will not say much on this front.

\section{Machine learning-based algorithms for problems in computational science}

In this and the next section, we will discuss how neural network models can be used to develop new algorithms.
For readers who are not familiar with neural networks, just think of them as being some replacement of polynomials.
We will discuss neural networks afterwards.

\subsection{Nonlinear multi-grid method and protein folding}

In traditional multi-grid method \cite{Brandt}, say for solving the linear systems of equation that arise from some
finite difference or finite element discretization of a linear elliptic equation, our objective is to  minimize a 
quadratic function like
\[ I_h(\bu_h) = \frac 12 \bu_h^T  L_h \bu_h - \bm{ f}_h^T \bu_h 
\]
Here $h$ is the grid size of the discretization.
The basic idea of the multi-grid method is to iterate between solving this problem and a reduced problem on a coarser grid with
grid size $H$.  In order to do this, we need the following 
\bi
\item a projection operator:  $P:\bu_h \rightarrow  \bu_H$, that maps functions defined on the fine grid to functions defined on
the coarse grid.
\item the effective operator at scale $H$:  $L_H = P^T L_h P$. This defines the objective function on the coarse grid:
\[ I_H(\bu_H) = \frac 12 \bu_H^T  L_H \bu_H - \bm{ f}_H^T \bu_H 
\]
\item a prolongation operator $Q: \bu_H \rightarrow  \bu_h$, that maps functions defined on the coarse grid to  functions defined on
the fine grid. Usually one can take $Q$ to be $P^T$.
\ei

The key idea here is coarse graining, and iterating between the fine scale  and the coarse-grained problem.
The main components in coarse graining is a set of coarse-grained variables and the effective coarse-grained problem.
Formulated this way, these are  obviously  general ideas that can be relevant for a wide variety of problems.
In practice, however, the difficulty lies in how to obtain the effective coarse-grained problem, a step that is trivial for
linear problems, and this is where
machine learning can help.

We are going to use the protein folding problem as an example to illustrate the general idea for nonlinear problems.

Let $\{\bx_j\}$ be the positions of  the atoms in a protein and the surrounding solvent, and 
$U=U(\{\bx_j \})$ be the potential energy of the combined protein-solvent system.
The potential energy consists of the energies due to chemical bonding, Van der Waals interaction, electro-static interaction, etc.
The protein folding problem is to find the ``ground state'' of the energy $U$:
 $$\mbox{``Minimize''} \, \, \, U. $$
 Here we have added quotation marks since really  what we want to do is to sample the distribution
 \[
   \rho_\beta= \frac1Z e^{- \beta U} , \, \beta = (k_B T)^{-1}
 \]

To define the coarse-grained problem, we assume that we are given a set of collective variables:
$\bm{s}=(s_1, \cdots, s_n), s_j = s_j(\bx_1,  \cdots, \bx_N),  (n <  N) $.
One possibility is to use the dihedral angles as the coarse-grained variables.
In principle, one may also use machine learning methods to learn the ``best'' set of coarse-grained variables
but this direction will not be pursued here.

Having defined the coarse-grained variables, the effective coarse-grained problem is simply the free energy associated with
this set of coarse-grained variables:
\[
  A(\bm s) = -\frac{1}\beta
  \ln p(\bm s), \quad
  p_\beta (\bm s)=
  \frac 1{Z} \int e^{-\beta U (\bm x)} \delta (\bm s(\bm x) - \bm s)\, d\bm x,
\]
Unlike the case for linear problems, for which the effective coarse-grained model is readily available, in the current situation,
 we have to find the function $A$ first.

The idea is to approximate $A$  by neural networks.  The issue here is how to obtain the training data.

Contrary to most standard machine learning problems where the training data is collected beforehand, 
in  applications to computational science and scientific computing, the training data is collected ``on-the-fly''
as learning proceeds.  This is referred to as the ``concurrent learning'' protocol \cite{EHanZhang2020}. In this regard, the standard machine 
learning problems for which the training data is collected beforehand are examples of  ``sequential learning''.
The key issue for concurrent learning is an efficient algorithm for generating the data in the best way.
The training dataset should on one hand be representative enough and on the other hand be as small as possible.

A general procedure for generating such datasets is suggested in \cite{EHanZhang2020}. 
It is called the EELT (exploration- examination-labeling-training) algorithm and it consists of the following steps:
     \begin{itemize} 
     \item exploration: exploring the  $\bm{s}$ space. This can be done 
      by sampling $\frac 1Z e^{-\beta A(\bm{s})}$ with the current approximation of $A$.
     \item examination: for each state explored, decide whether that state should be labeled. One way to do this is to use an  {\it a posteriori error estimator}.
     One possible such a posteriori error estimator is the variance of the predictions of an ensemble of machine learning models, see \cite{ZhangWangE2018}.
     \item labeling:  compute the mean force (say using restrained molecular dynamics)
    \[
    \bm F (\bm s) = -\nabla_{\bm s} A(\bm s).
    \]
    from which the free energy $A$ can be computed using standard thermodynamic integration.
      \item training:  train the appropriate neural network model.  To come up with a good neural network model, one has to take
      into account the symmetries in the problem.  For example, if we coarse grain a full atom representation of a collection of water
      molecules by eliminating the positions of the hydrogen atoms, then the free energy function for the resulting system should have
      permutation symmetry and this should be taken into account when designing the neural network model (see the next subsection).
      \end{itemize}

This can also be viewed as  a nonlinear multi-grid algorithm in the sense that it iterates between sampling $p_\beta$ on the space of the
coarse-grained variables
and the (constrained) Gibbs distribution $\rho_\beta$ for the full atom description.

This is a general procedure that should work for a large class of nonlinear ``multi-grid'' problems.

\begin{figure}
\centering
    \includegraphics[width=12cm]{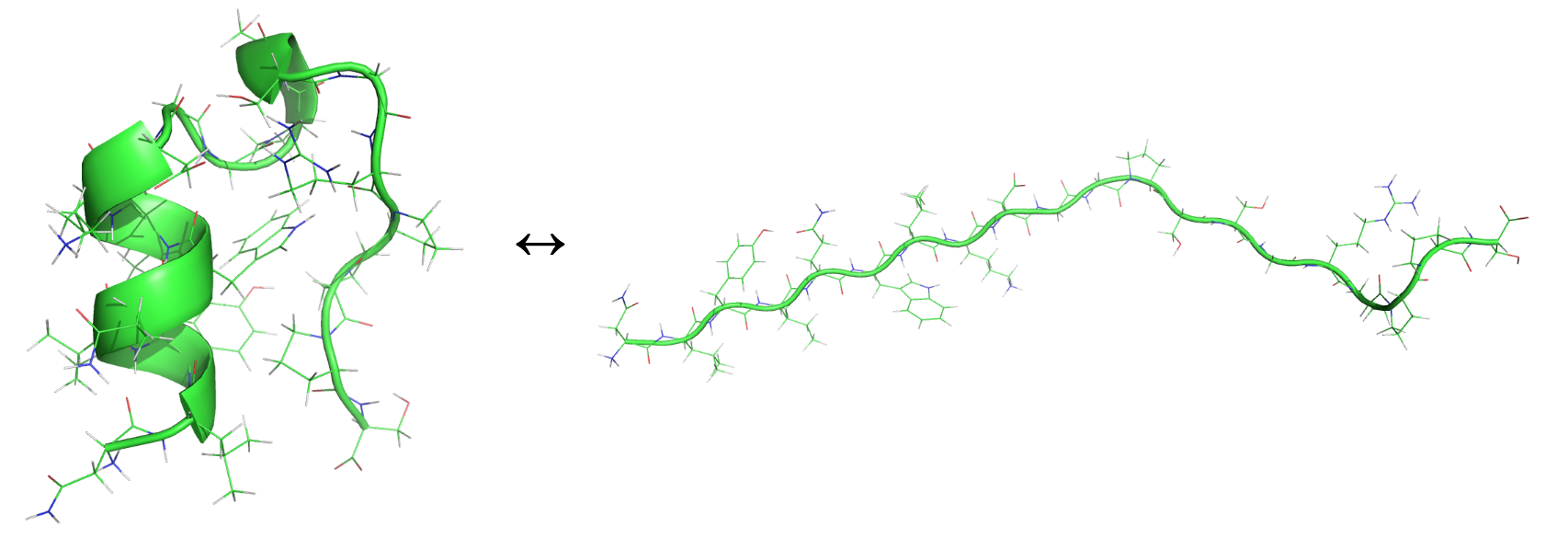}
   \caption{The folded and extended states of Trp-cage, reproduced with permission from \cite{WangZhangE2018}}.
     \label{fig:trp-cage}
\end{figure}

 Shown in Figure \ref{fig:trp-cage} is the extended and folded structure of  Trp-cage.
 This is a small protein with 20 amino acids. We have chosen the  38 dihedral angles as the collective variables.
 The full result is presented in \cite{WangZhangE2018}.
 
\subsection{Molecular dynamics with {\it ab initio} accuracy}

Molecular dynamics is a way of studying the behavior of molecular and material systems by tracking the trajectories of all the nuclei
in the system. The dynamics of the nuclei is assumed to obey Newton's law, with some  potential energy function (typically
called potential energy surface or PES) $V$ that models the effective interaction between the nuclei:
\[
m_i\frac{d^2\bm{x}_{i}}{dt^2}=-\nabla_{\bm{x}_{i}}V, \quad  V=V(\bm{x}_{1},\bm{x}_{2},...,\bm{x}_{i},...,\bm{x}_{N}),
\]

 How can we get the function $V$? Traditionally, there have been two rather different approaches.
The first is to compute  the inter-atomic forces ($-\nabla V$) on the fly using quantum mechanics models,
the most popular one being the density functional theory (DFT).  This is known as  the  Car-Parrinello  molecular
dynamics or {\it ab initio} molecular dynamics \cite{CPMD,E2011}. This approach  is quite 
accurate but is also very expensive, limiting the size of the system that one can handle to  about $1000$ atoms, even with
high performance supercomputers.
The other approach is to come up with empirical potentials.  Basically one guesses a functional form of $V$ with a small set of fitting parameters
which are then determined by a small set of data.
This approach is very  efficient but unreliable. 
This dilemma between accuracy and efficiency has been an essential road block for molecular dynamics for a long time.

With machine learning, the new paradigm is to use DFT to generate the data, and then use machine learning to generate an approximation
to  $V$.  This approach has the potential to produce an approximation to $V$ that is as accurate as the DFT model and as
efficient as the empirical potentials.

To achieve this goal, we have to address two issues. The first is the generation of data. The second is coming up with the
appropriate neural network model.  These two issues are the common features for all the problems that we discuss here.

The issue of adaptive data generation is very much the same as before.  The EELT procedure can still be used.
The details of how each step is implemented is a little different.  We refer to \cite{Zhangactive} for details.

\begin{figure}
\centering
    \includegraphics[width=15cm]{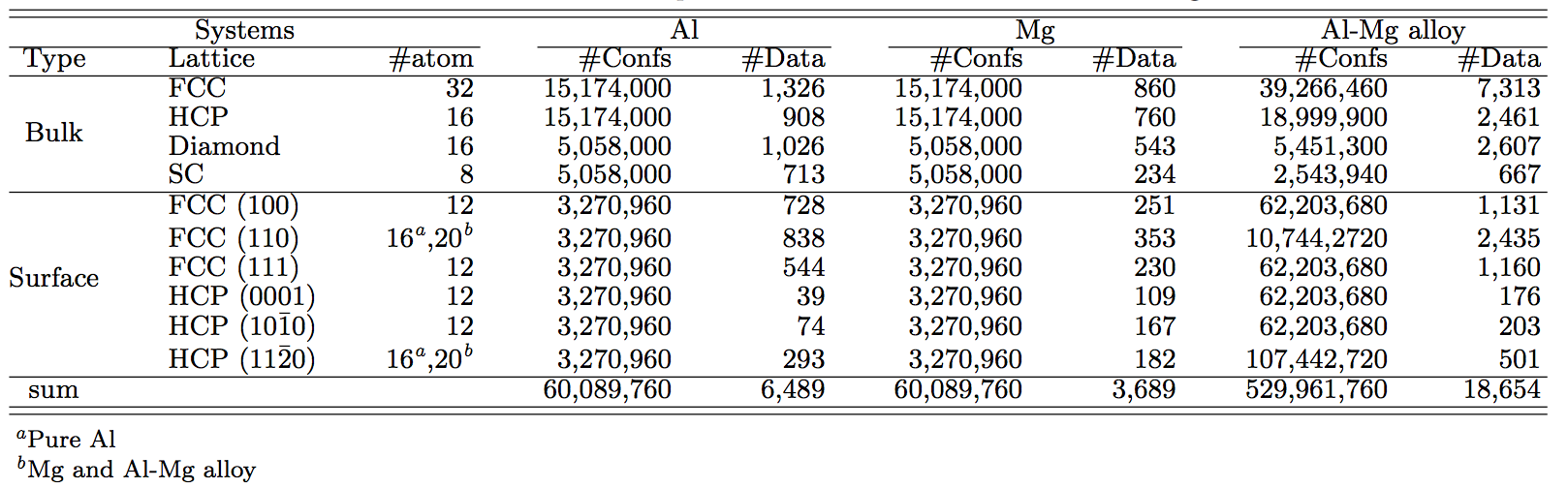}
    \caption{The results of EELT algorithm: Number of configurations explored vs. the number of data points labeled.
    Only a very small percentage of the configurations are labeled.  Reproduced with permission from Linfeng Zhang.
    See also \cite{Zhangactive}.}
       \label{fig:active-MD}
\end{figure}
Figure \ref{fig:active-MD} shows the effect of using the EELT algorithm.  As one can see, a very small percentage of the configurations
explored are actually labeled. For the Al-Mg example, only
$\sim$0.005\% configurations explored by  are selected for labeling.

\begin{figure}
\centering
 \includegraphics[width=15cm]{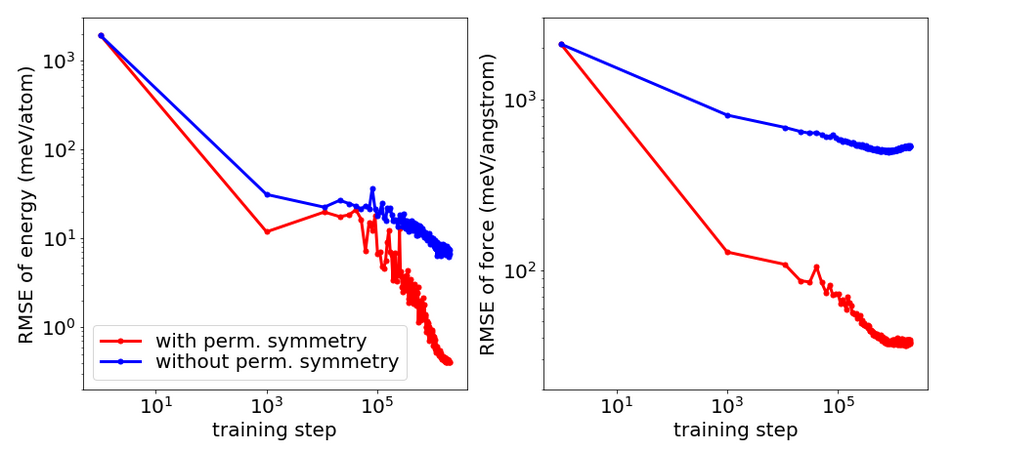}
 \caption{The effect of symmetry preservation on testing accuracy.  Shown in red are the results of poor man's way of imposing 
 symmetry (see main text for explanation). One can see that testing accuracy is drastically improved. Reproduced with permission
 from Linfeng Zhang. }
  \label{symm}

\end{figure}

For the second issue, the design of appropriate neural networks, the most important considerations are:
\begin{enumerate}
\item Extensiveness, the neural network should be extensive in the sense that if we want to extend the system, we just 
have to extend the neural network accordingly. One way of achieving this is suggested by Behler and Parrinello \cite{BP2007}.
\item Preserving the symmetry. Besides the usual translational and rotational symmetry, one also has the permutational
symmetry: If we relabel a system of copper atoms, its potential energy should not change.
It makes a big difference in terms of the accuracy of the neural network model whether one takes these symmetries into account
(see \cite{HanZhangCarE2017} and Figure \ref{symm}).
\end{enumerate}

One very nice and general way of addressing the symmetry problem is to design the neural network model as the composition of
two networks: An embedding network followed by a fitting network. The task for the embedding network is to represent
enough symmetry-preserving functions to be fed into the fitting network \cite{ZhangHanNIPS2018}.

\begin{figure}
\centering
  \includegraphics[width=15cm]{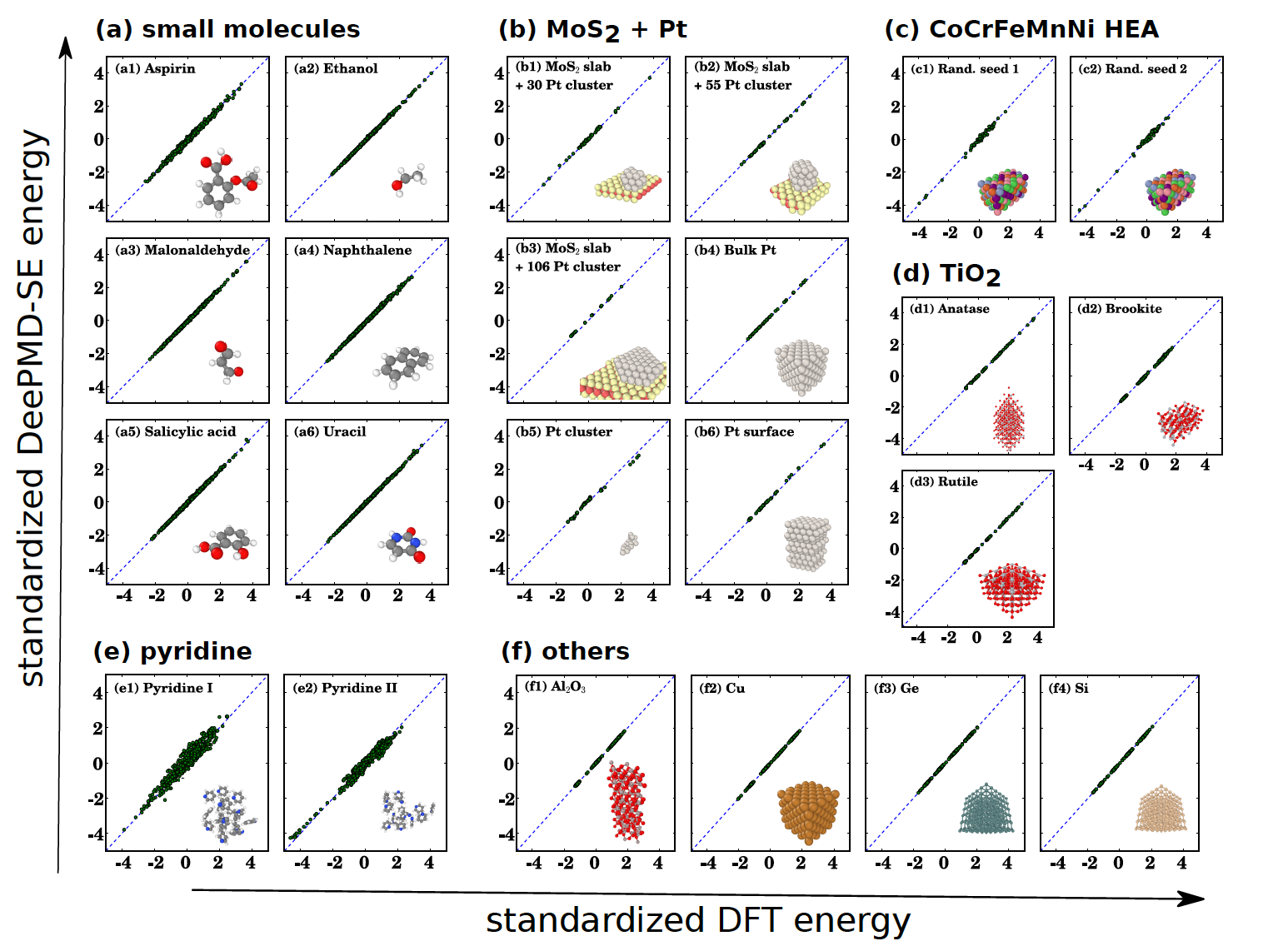}
\caption{The test accuracy of the Deep Potential for a wide variety of systems. Reproduced with permission from Linfeng Zhang.
  See also \cite{ZhangHanNIPS2018}.}
    \label{DeepP}
\end{figure}

With these issues properly addressed, one can come up with very satisfactory neural network-based representation of $V$
(see Figure \ref{DeepP}). This representation is named Deep Potential \cite{HanZhangCarE2017, ZhangHanNIPS2018}
and the Deep Potential-based molecular dynamics is named DeePMD \cite{DeePMD}.
As has been demonstrated recently in \cite{LuSC, JiaGB}, DeePMD, combined with state of the art high performance supercomputers,
can help to increase the size of the system that one can model with {\it ab initio} accuracy by  5 orders of magnitude.

\section{Machine learning-based algorithms for high dimensional problems in scientific computing}

\subsection{Stochastic control}

The first application of machine learning for solving high dimensional  problems in scientific computing  was presented in 
\cite{HanE2016}. Stochastic control was chosen as the first example due to its close analogy with machine learning.
Consider the discrete stochastic dynamical system:
\begin{equation}
  s_{t +1} = s_t+b_t(s_t,a_t)+\xi_{t +1}.
  \label{stoch-dyn}
	\end{equation}
Here	$s_t$ and $a_t$ are respectively the state and control at time $t$, $\xi_t$ is the noise at time $t$. Our objective is to solve:
\begin{align}
  \min_{\{a_t\}_{t =0}^{T-1}}\mathbb{E}_{\{\xi_t\}} \big\{\sum_{t =0}^{T-1} c_t(s_t,a_t)+c_T(s_T)\big\}
  \label{stoch-opt}
\end{align}
within the set of  feedback controls:
\begin{equation}
	a_t = A_t(s_t).
\end{equation}

We approximate the functions $A_t$ by neural network models:
	\begin{equation}
  	A_t(s)\approx \tilde{A}_t(s |\theta_t), \, t=0, \cdots, T-1
	\end{equation}
The optimization problem \eqref{stoch-opt} then becomes:
	\begin{equation}
  \min_{\{\theta_t\}_{t =0}^{T-1}}\mathbb{E}_{\{\xi_t\}} \big\{\sum_{t =0}^{T-1}c_t(s_t,\tilde{A}_t(s_t|\theta_t))+c_T(s_T)\}.
  \label{stoch-opt-NN}
	\end{equation}
	
\begin{figure}[H]
	\centering
	\includegraphics[width=0.9\textwidth]{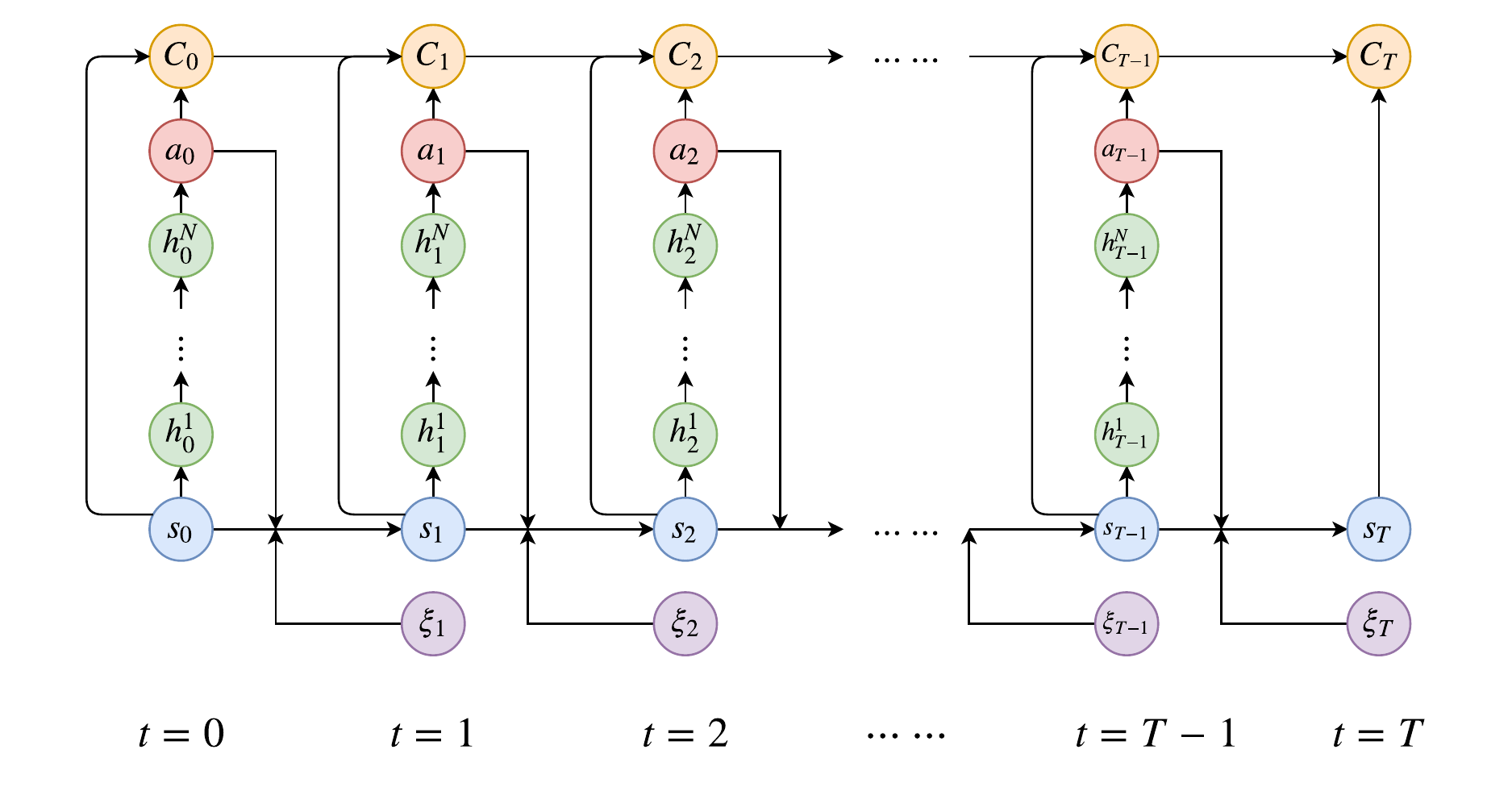}
	\vspace{-1em}
	\caption{Network architecture for solving stochastic control in discrete time. The whole network has $(N+1)T$ layers in total that involve free parameters to be optimized simultaneously. Each column (except $\xi_t$) corresponds to a sub-network at $t$. Reproduced with permission from
	Jiequn Han. See also \cite{HanE2016}.}
	\label{stoch-net}
	\end{figure}

	Unlike the situation in standard supervised learning, here we have $T$ set of neural networks to be trained simultaneously. 
	The network architecture is shown in Figure \ref{stoch-net}

Compared with the standard setting for machine learning, one can see a clear analogy in which 
\eqref{stoch-dyn} plays the role for the residual networks and the noise $\{\xi_t\}$ plays the role of data.
Indeed, stochastic gradient descent (SGD) can be readily used to solve the optimization problem \eqref{stoch-opt-NN}.

An example of the application of this algorithm is shown in Figure \ref{storage} for the problem of energy storage
with multiple devices. Here $n$ is the number of devices. For more details, we refer to \cite{HanE2016}.

\begin{center}
  \includegraphics[width=0.56\textwidth]{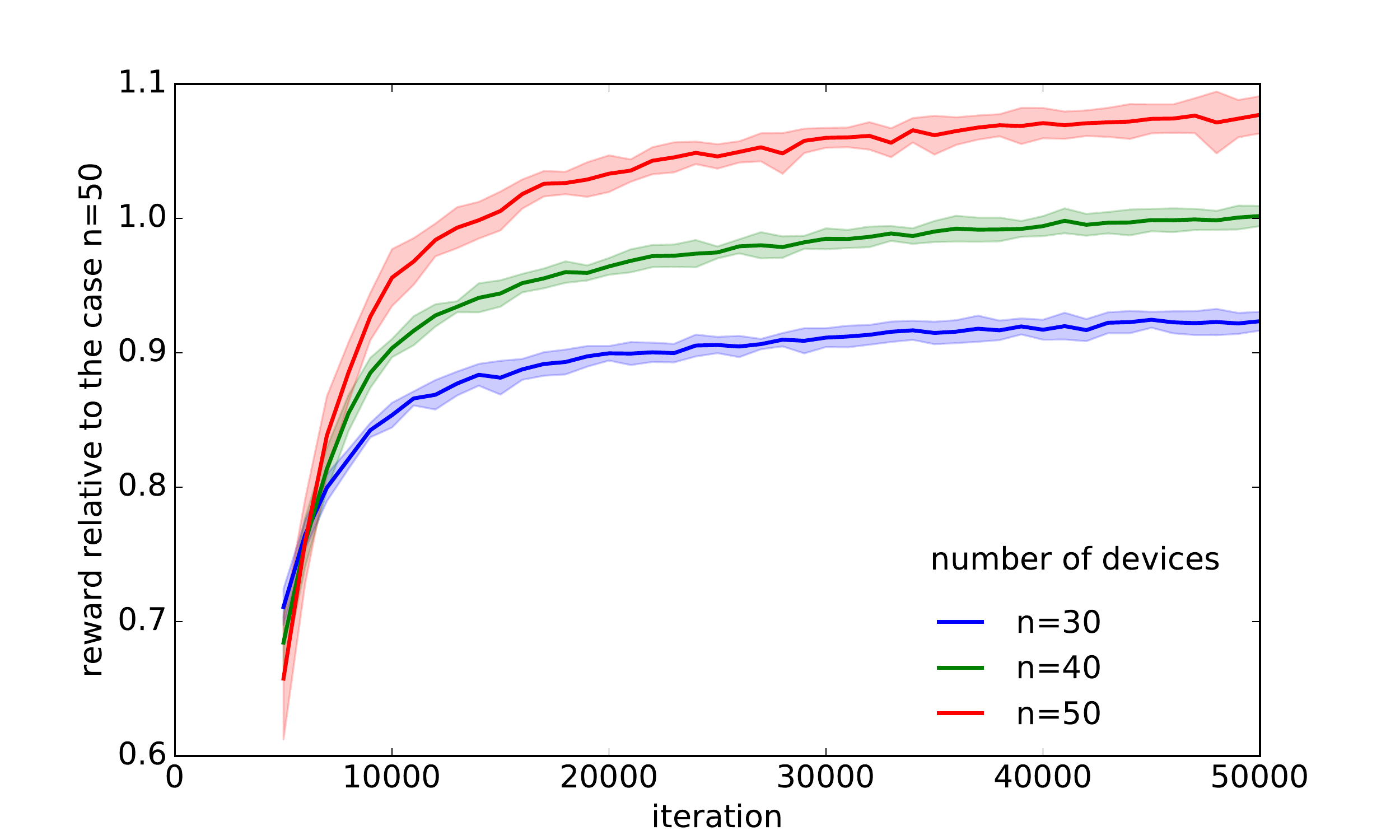}
  \captionof{figure}{Relative reward for the energy storage problem. The space of control function is $\mathbb{R}^{n+2}\rightarrow \mathbb{R}^{3n}$ for $n=30, 40, 50$, with multiple equality and inequality constrains. Reproduced with permission from Jiequn Han.  See also \cite{HanE2016}. }
  \label{storage}
\end{center}

\subsection{Nonlinear  parabolic PDEs}

Consider parabolic PDEs of the form:
\[
	  \frac{ \partial u}{ \partial t } 
	  + \frac{1}{2}  \sigma\sigma^T : \nabla^2_x u
	  +\mu \cdot \nabla u
	  +f\big(\sigma^{\operatorname{T}} \nabla u \big) = 0, \quad  u(T, x) = g(x)
\]
We study a terminal-value problem instead of initial-value problem since one of the main applications we have in mind is in finance.
To develop machine learning-based algorithms, we would like to first reformulate this  as a stochastic optimization problem.
This can be done using backward stochastic differential equations (BSDE) 
\cite{PardouxPeng1990}.
\begin{align}
&\inf_{Y_0,\{Z_t\}_{0\le t \le T}} \E |g(X_T) - Y_T|^2, \\
&s.t.\quad X_t = \xi + \int_{0}^{t}\mu(s,X_s)\, \,ds + \int_{0}^{t}\Sigma(s,X_s)\, dW_s, \\
&\hphantom{s.t.}\quad Y_t = Y_0 - \int_{0}^{t}h(s,X_s,Y_s,Z_s)\,  ds + \int_{0}^{t}(Z_s)^T dW_s.
\end{align}
It can be shown that the  unique minimizer of this problem is the solution to the PDE with:
\begin{equation}
\label{eq:nonlinear_Feynman_Kac}
  Y_t = u( t, X_t )
\qquad  
  \text{and}
\qquad 
  Z_t = \sigma^{T} ( t, X_t ) \, \nabla u( t, X_t ).
\end{equation}

With this formulation, one can develop a machine learning-based algorithm along the following lines,
adopting the ideas for the stochastic control problems discussed earlier
\cite{EHanJentzen2017,HanJentzenE2018}:
	\bi
		\item After time discretization, approximate the unknown functions
		$$X_0 \mapsto u(0, X_0)
		\qquad  
  		\text{and}
		\qquad 
		X_{t_j} \mapsto \sigma^{ \operatorname{T} }(t_j,X_{t_j}) \, \nabla u(t_j, X_{t_j})$$
		by feedforward neural networks $\psi$ and $\phi$.
		\item Using the BSDE, one constructs an approximation $\hat{u}$ that takes the paths 
		$\{ X_{ t_n } \}_{ 0 \leq n \leq N }$ and 
		$\{ W_{ t_n } \}_{ 0 \leq n \leq N }$ 
		as the input data and gives the final output,  
		denoted by 
		$
		  \hat{u}( 
		    \{ { X_{ t_n } } \}_{ 0 \leq n \leq N } , 
		    \{ W_{ t_n } \}_{ 0 \leq n \leq N } 
		  ) 
		$, 
		as an approximation to 
		$u( t_N, X_{ t_N } )$. 
		\item The error in the {matching between $\hat{u}$ and the given terminal condition} defines the expected loss function 
			\begin{equation*}
			  l(\theta) = 
			  \E\Big[
			    \big|g( X_{ t_N } ) - \hat{u}\big(\{ X_{ t_n } \}_{ 0 \leq n \leq N } , \{ W_{ t_n } \}_{ 0 \leq n \leq N }\big)\big|^2
			  \Big].
			\end{equation*}

			\ei
This algorithm is called the Deep BSDE method.

As applications, let us first study a stochastic control problem, but we now solve this problem using the Hamilton-Jacobi-Bellman (HJB) equation.
Consider the well-known LQG (linear quadratic Gaussian) problem at dimension $d=100$:
\begin{equation}
            dX_t = 2\sqrt{\lambda}\,m_t\,dt+\sqrt{2}\,dW_t,
        \end{equation}
 with the cost functional:
        $ J( \{ m_t \}_{ 0 \leq t \leq T } ) =
          \E\big[
            \int_0^T \|m_t\|_2^2 \, dt + g(X_T)
          \big]
        $.
The corresponding HJB equation is given by
        \begin{equation}
     \frac{ \partial u}{ \partial t }
  + \Delta u  - \lambda \|\nabla u \|_2^2 = 0
        \end{equation}
        Using the Hopf-Cole transform, we  can express the solution in the form:
\begin{equation}
                  u(t,x) = - \frac{ 1}{ \lambda }
                  \ln\!\bigg(
                    \E\Big[
                      \exp\!\Big(
                         - \lambda g( x + \sqrt{ 2 }W_{ T - t }  )
                      \Big)
                    \Big]
                  \bigg).
                \end{equation}
                This can be used to calibrate the accuracy of the Deep BSDE method.
                \begin{figure}[H]
                \vspace{-1em}
                \centering
	  \includegraphics[width=0.4\textwidth]{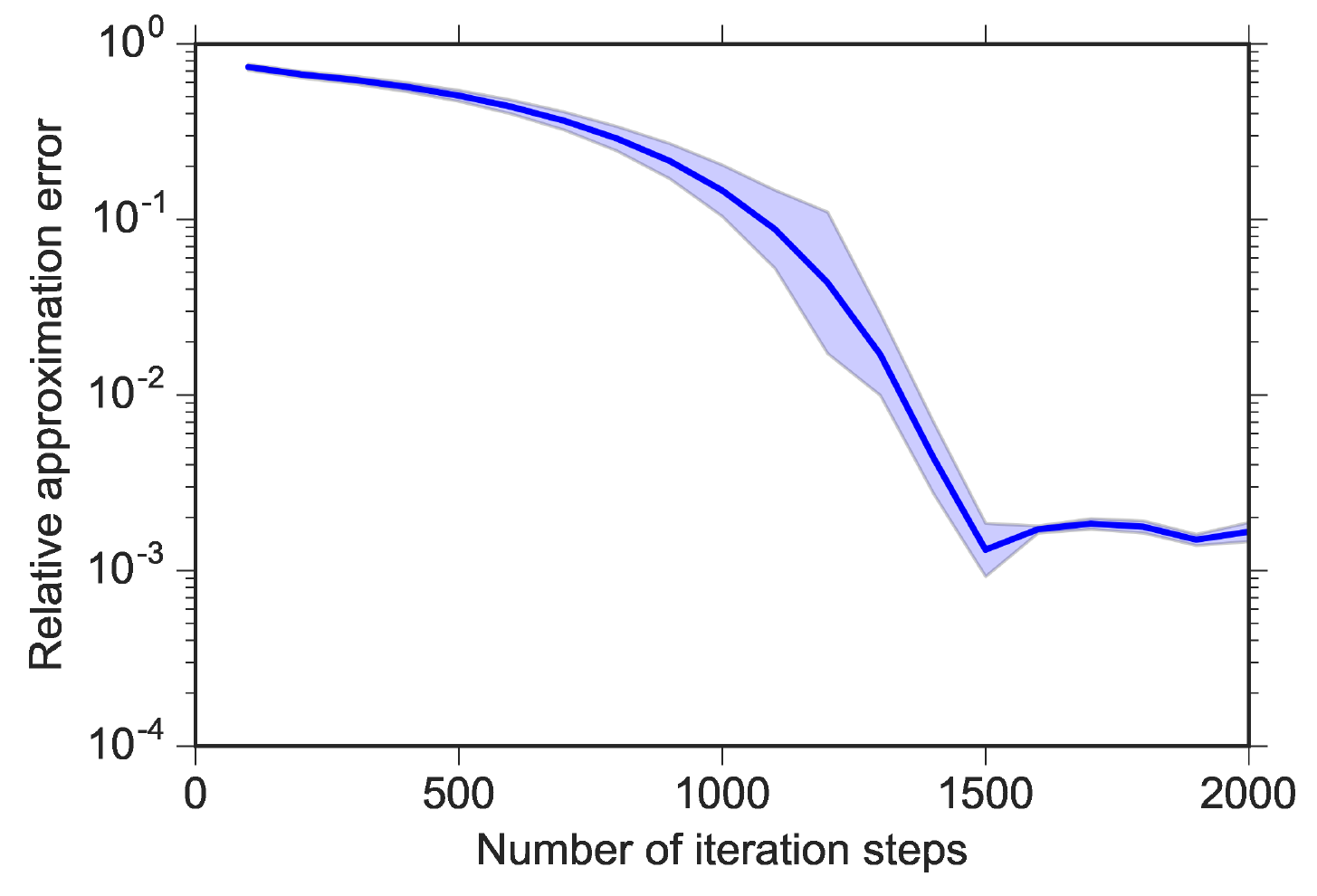}
	  \includegraphics[width=0.4\textwidth]{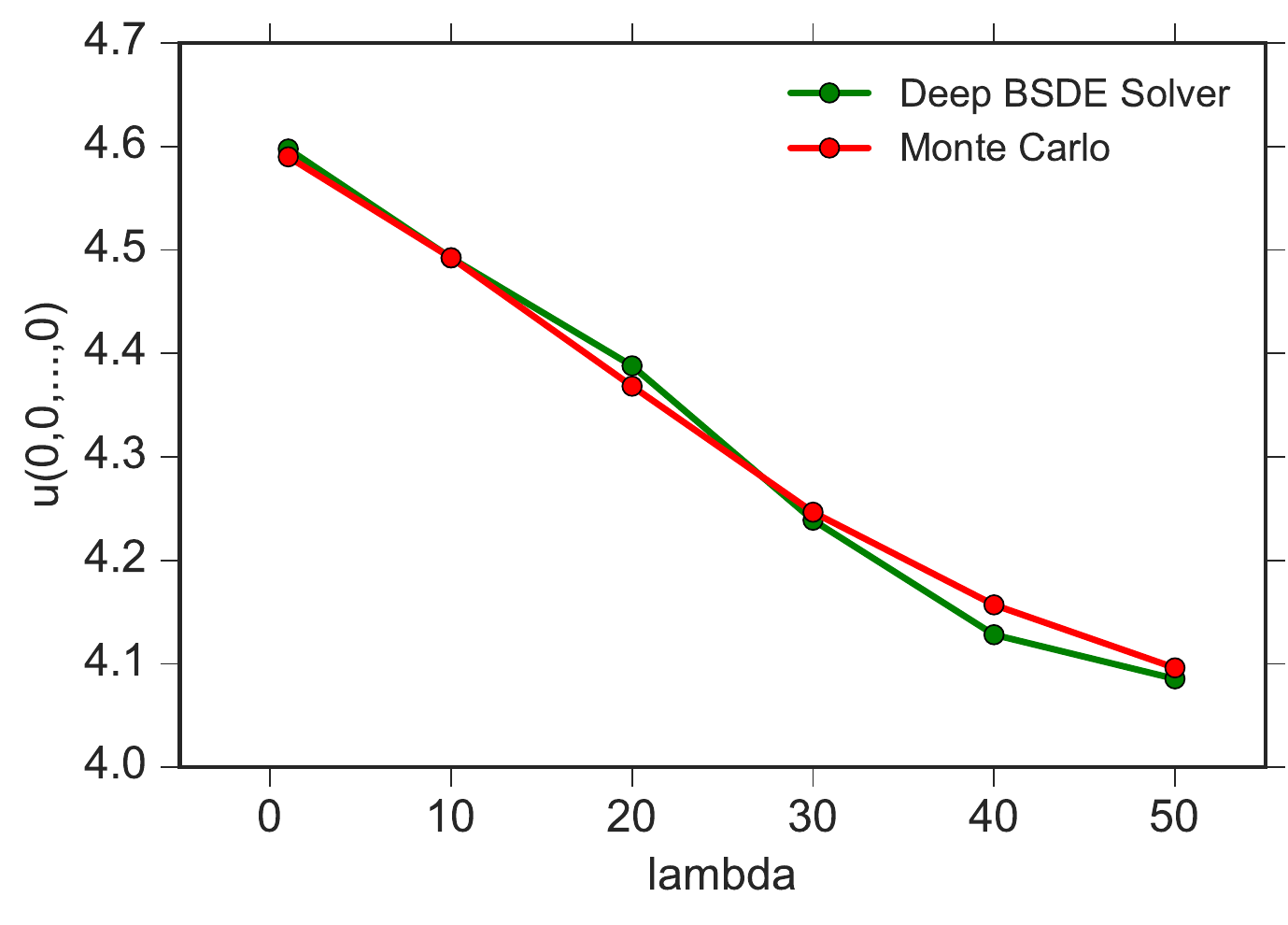}
	  \caption{{\footnotesize{Left: Relative error of the deep BSDE method for $ u( t{=}0, x{=}(0,\dots,0) )$ when $ \lambda = 1 $,
	  which achieves $ 0.17\% $ in a runtime of $ 330 $ seconds. Right: Optimal cost $u(t{=}0,x{=}(0,\dots,0))$ against different $\lambda$.}}
	  Reproduced with permission from Jiequn Han.  See also \cite{HanJentzenE2018}.}
		\end{figure}

As a second example, we study the Black-Scholes equation with default risk:
\[
		 \frac{\partial u}{\partial t} + \Delta u - \left( 1 - \delta \right) Q( u(t,x) ) \, u(t,x) - R \, u(t,x)= 0
		 \]
		 where $Q$ is some nonlinear function.
		 This form of modeling the default risk was proposed and/or used in the literature
		 for low dimensional situation ($d=5$, see \cite{HanJentzenE2018} for references).  
		 The Deep BSDE method was used for this problem with $d=100$ \cite{HanJentzenE2018}.
		 
		 \vspace{-.15in}
	\begin{figure}[H]
	\centering
	\includegraphics[width=0.6\textwidth]{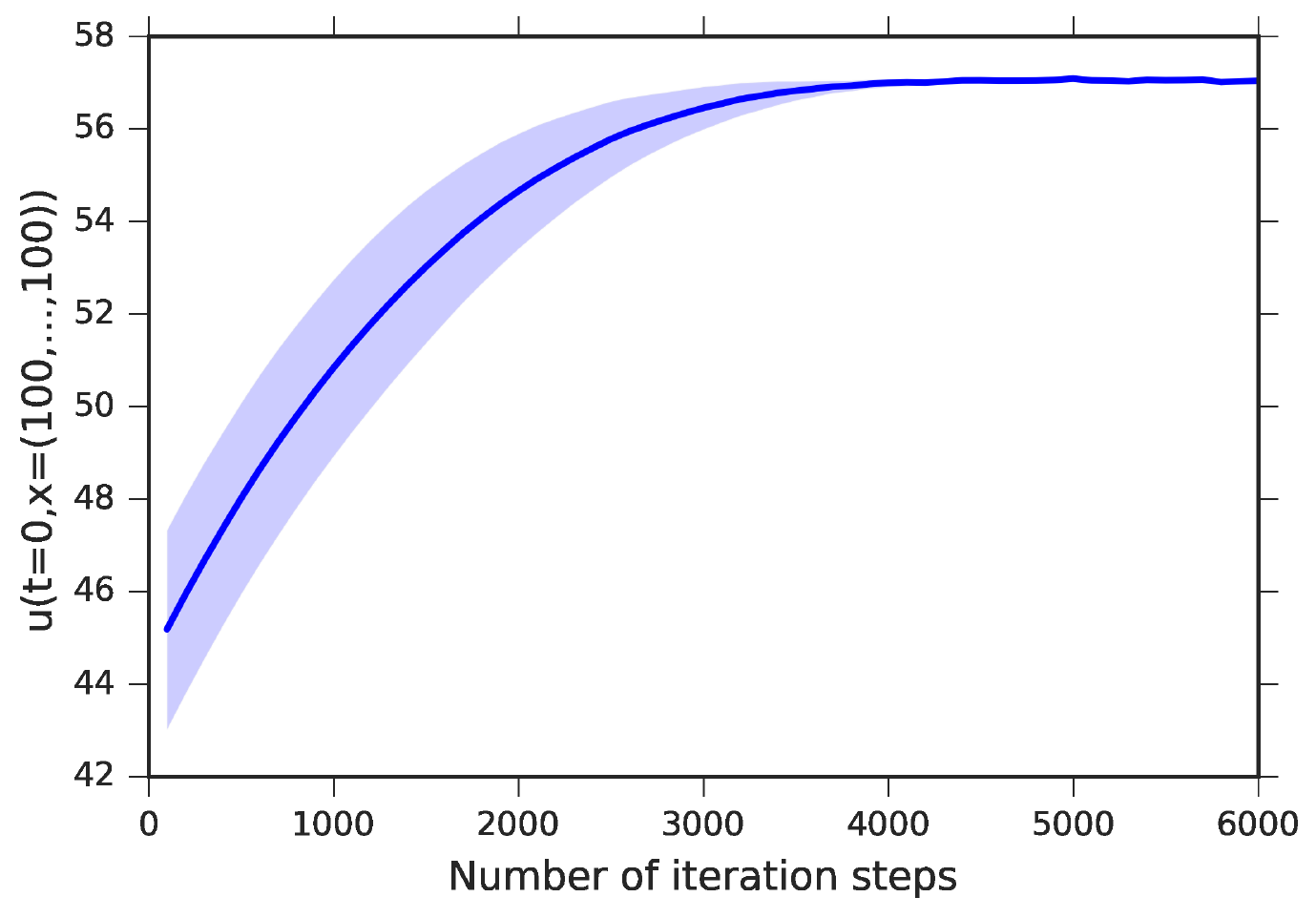}
	\caption{The solution of the Black-Scholes equation with default risk at  $d=100$.  
	The Deep BSDE method achieves a relative error of size $ 0.46\% $ in a runtime of $ 617 $ seconds. Reproduced with permission from Jiequn Han.  See also \cite{HanJentzenE2018}. }
	\end{figure}

The Deep BSDE method has been applied to  pricing basket options, interest rate-dependent options, Libor market model, Bermudan
Swaption, barrier option (see \cite{Review-1} for references).

\subsection{Moment closure for kinetic equations modeling gas dynamics}


The dynamics of gas can be modeled very accurately by the well-known 
{Boltzmann Equation}:
 \begin{equation}
  \dt f + \bv\cdot\dx f = \frac{1}{\veps} Q(f),  \quad \bv\in\R^3, \quad \bx\in\Omega\subset\R^3,
\end{equation}
where $f$ is the  phase space density function,  $\varepsilon$ is the Knudsen number:
       \begin{equation*}
        \varepsilon = \frac{\text{mean free path}}{\text{macroscopic length}},
      \end{equation*}
 $Q$ is  the collision operator
that models the collision process between gas particles.
When $\veps \ll 1$, this can be approximated by the Euler equation:
\begin{equation}
  \dt \bU + \dx\cdot \bF(\bU) = 0,
\end{equation}
where $$\bU = (\rho, \rho\bu, E)^T,  \quad
 \rho = \int f\diff\bv, \quad \bu = \frac{1}{\rho}\int f\bv\diff\bv.
 $$
and
$$\bF(\bU) = (\rho\bu, \rho \bu\otimes\bu + pI, (E + p)\bu)^T $$
Euler's equation can be obtained by projecting the Boltzmann equation on to the low order moments
involved, and making use of the ansatz that the distribution function $f$ is close to be local Maxwellian.

 \begin{figure}[H]
    \centering
    \includegraphics[width=\textwidth]{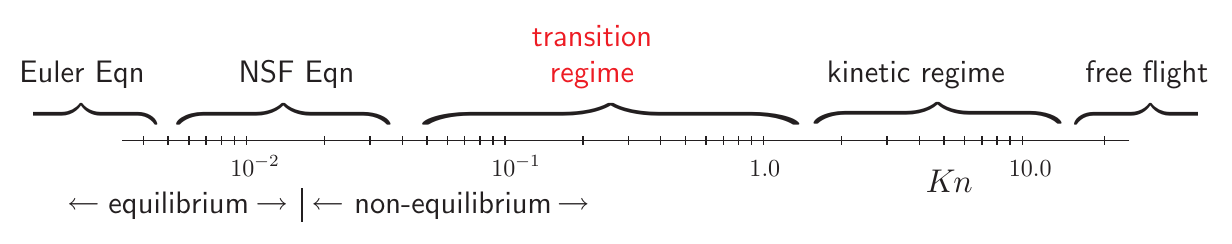}
    \caption{The different regimes of gas dynamics. Reproduced with permission from Jiequn Han.  See also \cite{HMME2019}.}
  \end{figure}

What happens when $\veps$ is not small?  In this case, a natural idea is to seek  
some generalization of Euler's equation using more moments.
This program was initiated by Harold Grad who constructed the well-known 
 13-moment system  using the moments of  $\{1, \bv, (\bv - \bu)\otimes(\bv - \bu), |\bv - \bu|^2(\bv - \bu)\}$.
 This line of work has encountered several difficulties. First, there is no guarantee that the
 equations obtained are well-posed.  Secondly there is always the ``closure problem'': When projecting the
 Boltzmann equation on a set of moments, there are always terms which involve moments outside the set of
 moments considered. In order to obtain a closed system, one needs to model these terms 
in some way.  For Euler's equation, this is done using the local Maxwellian approximation.
This is accurate when $\veps$ is small, but is no longer so when $\veps$ is not small.
It is highly unclear what should be used as the replacement.

In \cite{HMME2019}, Han et al developed a machine learning-based moment method.
The overall objective is to construct an uniformly accurate (generalized) moment model. The methodology consists of two steps:

1: Learn a set of optimal generalized moments through an auto-encoder.
Here by optimality we mean that the set of generalized moments retains a maximum amount of information about the original
distribution function and can be used to recover the distribution function with a minimum loss of accuracy.
This can be done as follows:
Find an encoder $\Psi$ and  a decoder $\Phi$ that recovers the original $f$ from $\bU, \bW$
\[
\bW=\Psi(f)=\int \bw f\diff \bv, \quad \Phi(\bU, \bW)(\bv) = \bm{h}(\bv; \bU, \bW).
\]
\[
\mbox{Minimize}_{\bw, \bm{h}}  \, \underset{f\sim \mathcal{D}}{\E} \|f - \Phi(\Psi(f))\|^2. 
\]
$\bU$ and $\bW$ form the  set of generalized hydrodynamic variables that we will use to model the gas flow.

{2: Learn the  fluxes and source terms in the PDE for the projected PDE.}
The effective PDE for $\bU$ and $\bW$ can be  obtained by formally projecting the Boltzmann equation on this set of (generalized) moments.
This gives us a set of PDEs of the form:
\begin{equation}
  \left\{
  \begin{aligned}
     & \dt \bU + \dx\cdot \bF(\bU, \bW; \veps) = 0, \\
     & \dt \bW + \dx\cdot \bG(\bU, \bW; \veps) = \bm{R}(\bU, \bW; \veps).
  \end{aligned}
  \right.
\end{equation}
where  $\bF(\bU, \bW; \veps)= \int \bv \bU f d\bv,  \bG(\bU, \bW; \veps)= \int \bv \bW f d \bv,
 \bm{R}(\bU, \bW; \veps) =  {\veps}^{-1} \int \bW Q(f) d \bv $.
Our task now is  to learn  $\bF, \bG, \bm{R}$ from the original kinetic equation.

Again the important issues are (1) get an optimal dataset, and (2) enforce the physical constraints.
Here two notable physical constraints are (1) conservation laws and (2) symmetries.  Conservation laws are
automatically respected in this approach.
Regarding symmetries, besides the usual static symmetry, there is now a new dynamic symmetry: the Galilean invariance.
These issues are all discussed in \cite{HMME2019}.
We also refer to \cite{HMME2019} for  numerical results for the models obtained this way.

\section{Mathematical theory of machine learning}

While neural network-based machine learning has demonstrated a wide range of very impressive successes,
it has also acquired a reputation of being a  ``black magic'' rather than a solid scientific technique.
This is due to the fact that (1) we lack a basic understanding of the fundamental reasons behind its success;
(2) the performance of these models and algorithms is quite sensitive to the choice of the hyper-parameters such 
as the architecture of the network and the learning rate; and (3) some techniques, such as batch normalization, does
appear to be a black magic.

To change this situation, we need to (1) improve our understanding of the reasons behind the success and the fragility of
neural network-based models and algorithms and (2) find ways to formulate more robust models and design more 
robust algorithms.  In this section we address the first issue.  The next section will be devoted to the second issue.

Here is a list of the most basic questions that we need to address:
\bi
\item Why does it work in such high dimension? 
\item Why simple gradient descent works for training neural network models? 
\item Is over-parametrization good or bad? 
\item Why does neural network modeling require such extensive parameter tuning?
\ei
At this point, we do not yet have clear answers to all these questions.  But some coherent picture is starting to emerge.
We will focus on the problem of supervised learning, namely approximating a target function using a finite dataset.
For simplicity, we will limit ourselves to the case when the physical domain of interest is $X=[0, 1]^d$.

\subsection{An introduction to neural network-based supervised learning}

The basic problem in supervised learning is as follows:
Given a natural number $ n \in \N $
and a sequence 
$ \{( \bx_j, y_j ) = ( \bx_j , f^*( \bx_j ) ) $,
$ j \in \{ 1, 2, \dots, n \} \}$, 
of pairs of input-output data,  
we want to  recover the target function $ f^* $ as accurately as possible. 
We will assume that the input data $\{ \bx_j,  j \in \{ 1, 2, \dots, n \} \}$, 
is sampled from the probability distribution $ \mu $ on $ \R^d $.

{\bf Step 1. Choose a hypothesis space}. This is a set of trial functions $ \cH_m $
where $ m \in \N $ is  the dimensionality of $\cH_m$.
One might choose piecewise polynomials or wavelets.
In modern machine learning the most popular  choice is neural network functions.
Two-layer neural network functions (one input layer, one output layer which usually does not count, and one hidden layer) take the form 
\begin{equation}
f_m(\bx,  \theta) =
\frac 1 m \sum_{ j = 1 }^m a_j \sigma( \left< \bw_j , \bx \right> )
\end{equation}
where $ \sigma \colon \R \to \R $ is a fixed scalar nonlinear function 
and where $ \theta =\{ ( a_j, w_j )_{ j \in \{ 1, 2, \dots, m \} } \}$ are the parameters to be optimized (or trained). 
A popular example for the nonlinear function $ \sigma \colon \R \to \R $ is the  
 ReLU (rectified linear unit) activation function:
$ \sigma(z) = \max\{ z, 0 \}$, for all $ z \in \R $.  We will restrict our attention to this activation function.
Roughly speaking, deep neural network (DNN) functions are obtained if one composes two-layer neural network functions several times.
One important class of DNN models are residual neural networks (ResNet). They closely resemble discretized 
ordinary differential equations and take the form 
\begin{align}
\label{ResNet}
\bz_{l+1} & = \bz_l +  \sum_{j=1}^M \ba_{j,l}\sigma( \left<\bz_l, \bw_{j,l} \right> ), 
\qquad
\bz_0 = \bV \tilde{\bx}, 
\qquad
f_L(\bx, \theta) = \left< \alpha, \bz_L \right>
\end{align}
for $ l \in \{ 0, 1, \dots, L - 1 \} $ where $ L, M \in \N $. 
Here the parameters are 
$ 
  \theta 
  = ( \alpha, \bV, 
    ( \ba_{j,l} )_{ j, l } 
    , 
    ( \bw_{ j, l } )_{ j, l }
  ) 
$.
ResNets are the model of choice for truly deep neural network models. 

{\bf Step 2. Choose a loss function}. The primary consideration for the
choice of the loss function is to fit the data.  Therefore the most obvious
choice is the $L^2$ loss:
\begin{equation}
\hat{\mathcal{R}}_n(f) = \frac1n \sum_{j=1}^n | f(\bx_j) - y_j |^2 = \frac 1n \sum_{j=1}^n | f(\bx_j) - f^*(\bx_j) |^2.
\label{empirical-risk}
\end{equation}
This is also called the ``empirical risk''.
Sometimes we also add regularization terms.

{\bf Step 3. Choose an optimization algorithm}.  
The most popular optimization algorithms in machine learning are different versions of
the gradient descent (GD) algorithm, or its stochastic analog, the
stochastic gradient descent (SGD) algorithm. Assume that the objective function we aim to minimize is of the form
\begin{equation}
F(\theta) = \E_{\xi \sim \nu}\big[ l(\theta, \xi) \big].
\label{expectation}
\end{equation}
The simplest form of SGD iteration takes the form 
\begin{equation}
\theta_{k+1} = \theta_k - \eta \nabla l(\theta_k, \xi_k),
\end{equation}
for $ k \in \N_0 $ where $\{ \xi_k $, $ k \in \N_0 = \{ 0, 1, 2, \dots \} \}$ is a sequence of 
i.i.d.\ random variables sampled from the distribution $ \nu $ 
and $ \eta $ is the learning rate which might also change during the course of the iteration.
In contrast, GD takes the form
\begin{equation}
\theta_{k+1} = \theta_k - \eta \nabla \E_{\xi \sim \nu}\big[ l(\theta_k, \xi) \big].
\end{equation}
Obviously this form of SGD can be adapted to loss functions of the form \eqref{empirical-risk} which can  be
regarded as an expectation with $\nu$ being the empirical distribution on the training dataset.
This DNN-SGD paradigm is the heart of modern machine learning.

\subsection{Approximation theory}

The simplest way of approximating functions is to use polynomials.
For polynomial approximation, there are two kinds of theorems.
The first is the Weierstrass' Theorem which asserts that  continuous functions can be uniformly approximated by polynomials
on compact domains. 
The second is  Taylor's Theorem which tells us that the rate of convergence  depends on the smoothness of the target function.

Using the terminology in neural network theory, Weierstrass' Theorem is the ``Universal Approximation Theorem'' (UAT).
It is a useful fact. 
But Taylor's Theorem is more useful since it tells us something about the rate of convergence.
The form of Taylor's Theorem used in approximation theory 
 are the direct and inverse approximation theorems which assert that a given function can approximated
by polynomials with a particular rate if and only if certain norms of that function is finite.
This particular norm, which measures the regularity of the function, is the key quantity that characterizes this approximation scheme.
For piecewise polynomials, these norms are some Besov space norms \cite{Triebel}.  For $L^2$, a typical result looks like follows:
\beq
\inf_{f \in \mathcal{H}_m} \|f - f_m\|_{L^2(X)} \le C_0 h^{\alpha} \|f\|_{H^{\alpha}(X)}
\label{rate-CoD}
\eeq
Here $H^{\alpha}$ stands for the $\alpha$-th order Sobolev norm \cite{Triebel}, $m$ is the number of degrees of freedom.
On a regular grid, the grid size is given by  
\beq
h \sim m^{-1/d}
\eeq

An important thing to notice is that the convergence rate in \eqref{rate-CoD}  suffers from CoD:  If we want to reduce the
error by a factor of $\epsilon$, we need to increase $m$ by a factor $m \sim \epsilon^{-d}$ if $\alpha=1$.
For $d=100$ which is not very high dimension by the standards of machine learning, this means that we have to increase
$m$ by a factor of $\epsilon^{-100}$. This is why polynomials and piecewise polynomials are not useful in high dimensions.

Another way to appreciate this is as follows. The number of monomials of degree $p$ in dimension $d$ is $C^d_{p+d}$.
This grows very fast for large values of $d$ and $p$.

{What should we expect in high dimension?}
One example that we can learn from is Monte Carlo methods for integration. Consider the problem of approximating 
$$I(g) = \EE_{\bx \sim \mu} g(\bx) $$
using 
$$ I_m(g) = \frac 1m  \sum_j g(\bx_j) $$
where $ \{\bx_j, j  \in [m] \}$ is a set of  i.i.d samples of  the probability distribution $\mu$. A direct computation gives
$$\EE(I(g)-I_m(g))^2 = \frac{\mbox{var}(g)}m, \quad 
\mbox{var}(g) =  \EE_{\bx \sim \mu} g^2(\bx)  - ( \EE_{\bx \sim \mu} g(\bx))^2 
$$
This exact relation tells us two things. (1) The convergence rate of Monte Carlo integration is independent of dimension.
(2) The error constant is given by the variance of the integrand. Therefore to reduce error, one has to do variance reduction.

Had we used grid-based quadrature rules, the accuracy would have also suffered from CoD.

It is possible to improve the Monte Carlo rate by more sophisticated ways of choosing $ \{\bx_j, j  \in [m] \}$, say using
number-theoretic-based quadrature rules. But these typically give rise to an $O(1/d)$ improvement for the convergence rate
and it diminishes as $d \rightarrow \infty$.

Based on these considerations, we aim to find function approximations that satisfy:
\[
 \inf_{f \in \mathcal{H}_m}  \mathcal{R}(f)  = 
\inf_{f \in \mathcal{H}_m}  \|f - f^*\|^2_{L^2(d \mu)}
\lesssim \frac{\|f^*\|_*^2}{m} 
\]
The natural questions are then:
\bi
\item How can we achieve this? That is, what kind of hypothesis space should we choose?
\item What should be the ``norm'' $\|\cdot \|_* $ (associated with the choice of $\mathcal{H}_m$)?
Here we put norm in quotation marks since it does not have to be a real norm. All we need is that it controls the
approximation error.
\ei

{Regarding the first question, let us look an illustrative example.}

Consider the following Fourier representation of the function $f$ and its approximation $f_m$ (say FFT-based):
\begin{equation*}
f(\bx) = \int_{\R^d} a(\bom) e^{i (\bom, \bx)} d \bom,
\quad
f_m(\bx) =  \frac 1m \sum_j a(\bom_j) e^{i (\bom_j, \bx)} 
\end{equation*}
Here $\{\bom_j\}$ is a fixed grid, e.g. uniform grid.  For this approximation, we have
$$
 \|f - f_m\|_{L^2(X)} \le C_0 m^{-\alpha/d} \|f\|_{H^{\alpha}(X)}
$$
which suffers from CoD.

Now consider the alternative representation 
\begin{equation}
f(\bx) = \int_{\R^d} a(\bom) e^{i (\bom, \bx)}  \pi(d\bom) =
\E_{\bom \sim \pi} a(\bom) e^{i (\bom, \bx)}
\label{Fourier-alternative}
\end{equation}
where $\pi$ is a probability distribution.  Now to approximate $f$, it is natural to use Monte Carlo.
Let $\{\bom_j\}$ be an i.i.d. sample of $\pi$,  $f_m(\bx) = \frac 1m \sum_{j=1}^m a(\bom_j) e^{i (\bom_j, \bx)}$,   then we have
$$
\E|f(\bx) - f_m(\bx)|^2 =\frac{ {\mbox{var}(f)} } m
$$
This approximation does not suffer from CoD. Notice that 
$ f_m(\bx) = \frac 1m \sum_{j=1}^m a_j \sigma (\bom_j^T \bx)$  is nothing but a   two-layer neural network  with activation function
 $\sigma(z) = e^{iz}$ (here $a_j = a(\bom_j)$).
 
 We believe that this simple argument is really at the heart of why neural network models do so well in high dimension.

Now let us turn to a concrete example of the kind of approximation theory for neural network models.
We will consider two-layer neural networks.
$$ \cH_m = \{ f_m(\bx) =
\frac 1 m \sum_j  a_j \sigma(\wb_j^T \xb) \}, \, \theta=\{(a_j, \wb_j), j \in [m] \}
$$
Consider  function $f: X = [0,1]^d\mapsto \RR$ of the following form
$$ f(\bx) = \int_{\Omega} a \sigma (\bw^T  \bx ) \rho (da, d\bw)
=\EE_{(a, \wb) \sim \rho}[a\sigma(\bw^T\bx)], \quad \bx \in X$$
where $\Omega = \RR^1 \times \RR^{d+1} $,
$\rho$ is a probability distribution on $\Omega$.
Define:
$$\|f \|_{\mathcal{B}} = \inf_{\rho\in P_f} \left(\EE_\rho [a^2\|\bw\|^2_1] \right)^{1/2} $$
where $P_f:=\{\rho: f(\bx)=\EE_\rho[a\sigma(\bw^T\bx)]\}$. This is called the Barron norm \cite{Bach2017,EMaWu2018, EMaWu2019} .
The space
$$
\mathcal{B} =  \{ f \in C^0 : \|f\|_{\mathcal{B}} < \infty \}
$$
is called the Barron space \cite{Bach2017,EMaWu2018, EMaWu2019} (see also \cite{Barron1993, KlusowskiBarron2016,  EWojtowytsch2020}).

In analogy with classical approximation theory, we can also prove some direct and inverse approximation theorem \cite{EMaWu2019}.

\begin{theorem}[Direct Approximation Theorem]
If $ \|f \|_{\mathcal{B}} < \infty$, then for any integer $m > 0$, there exists a two-layer neural network function $f_m$ such that
$$
\| f- f_m \|_{L^2(X)} \lesssim  \frac{ \|f \|_{\mathcal{B}}}{\sqrt{m}}
$$
\end{theorem}

\begin{theorem}[Inverse Approximation Theorem]
Let
\[
\mathcal{N}_{C} \Def \{\,\frac{1}{m}\sum_{k=1}^m a_k\sigma(\bw_k^T\bx) : 
\frac{1}{m}\sum_{k=1}^m |a_k|^2\|\bw_k\|_1^2 \leq C^2, m\in \mathbb{N}^{+}\,\}.
\]
Let $f^*$ be a continuous function. Assume
there exists a constant $C$ and a sequence of functions $f_m \in \mathcal{N}_{C} $
such that
$$ f_m(\bx) \rightarrow f^*(\bx) $$
for all $\bx \in X$. Then
there exists a probability distribution $\rho^*$ on $\Omega$, such that
\[
   f^*(\bx) = \int a\sigma(\bw^T\bx) \rho^*(da,d\bw),
\]
for all $\bx \in X$ and $\|f^*\|_{\mathcal{B}} \le C$.
\end{theorem}

\subsection{Estimation error}

\begin{figure}
\centering
      \includegraphics[width=10cm]{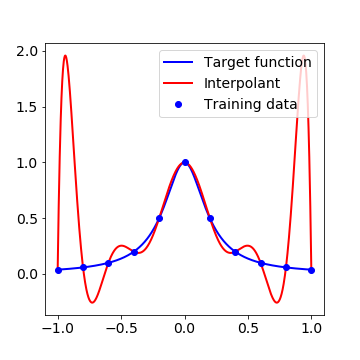}
    \caption{The Runge phenomenon: $f^*(x)= \frac {1}{1+ 25 x^2}$. Reproduced with permission from Chao Ma.}
    \label{fig:runge}
\end{figure}

Another issue we have to worry about is the performance of the machine learning model outside the training dataset.
This issue also shows up in classical approximation theory.  
Illustrated in Figure \ref{fig:runge} is the classical
Runge phenomenon for polynomial interpolation on a uniform grid.  One can see that away from the grid points,
the error of the interpolant can be very large. This is a situation that we would like to avoid.

What we do in practice is to minimize the training error:
$$\hat{\mathcal{R}}_n (\theta) = \frac 1n \sum_j (f(\bx_j, \theta) - f^*(\bx_j))^2 $$
but we are interested in the testing error, which is a sampled version of the population risk:
$${\mathcal{R}} (\theta) = \E_{\bx \sim \mu} (f(\bx, \theta) - f^*(\bx))^2 $$
The question is how we can control the difference between these two errors.

One way of doing this is to use the notion of Rademacher complexity.
The important fact for us here is that the Rademacher complexity controls the difference between training and testing
errors (also called the ``generalization gap''). Indeed, let
 $\mathcal{H}$ be  a set of functions,  and $S=(\bx_1, \bx_2, ..., \bx_n)$ be a dataset. Then,
up to logarithmic terms,  we have
    \[
      \sup_{h\in\mathcal{H}} \left| \E_\bx \left[h(\bx)\right] - \frac{1}{n}
      \sum_{i=1}^n h(\bx_{i}) \right|
      \sim  \rad_S(\mathcal{H}) 
          \]
where the {\bf Rademacher complexity} of $\mathcal{H}$ with respect to $S$ is defined as
\begin{equation}
\rad_S(\mathcal{H}) = \frac{1}{n}\mathbb{E}_\xi \left[\sup_{h\in\mathcal{H}}\sum\limits_{i=1}^n \xi_ih(\bx_i)\right],
\end{equation}
where $\{\xi_i\}_{i=1}^n$ are i.i.d. random variables taking values $\pm1$ with equal probability.

The question then becomes to bound the Rademacher complexity of a hypothesis space.  For the Barron space, we have \cite{Bach2017}:

\begin{theorem} 
Let $\cF_Q = \{ f \in \mathcal{B}, \|f\|_{\mathcal{B}} \le Q  \}$.
Then we have
\[
    \rad_S(\cF_Q) \leq 2Q \sqrt{\frac{2\ln(2d)}{n}}
\]
where $n = |S|$, the size of the dataset $S$.
\end{theorem}

\subsection{A priori estimates for regularized models}

{Consider the  regularized model} 
     \begin{equation}
\cL_n(\theta) = \hat{\mathcal{R}}_n(\theta) +
\lambda \sqrt{\frac{\log(2d)}{n}} \|\theta\|_{\mathcal{P}}, \quad \quad
        \hat{\theta}_n = \argmin \cL_n(\theta)
    \end{equation}
where the path norm is defined by:
    \[
        \|\theta\|_{\mathcal{P}} = \left(\frac 1m \sum_{k=1}^m |a_k|^2\|\bw_k\|_1^2 \right)^{1/2}
    \]

The following result was proved in \cite{EMaWu2018}:

  \begin{theorem}
Assume  $f^*: X \mapsto [0,1] \in \mathcal{B}$.
There exist constants  $C_0$, such that for any $\delta > 0$,
        if $\lambda\geq C_0$, then with probability  at least $1-\delta$ over the choice of the training dataset, we have
        \[
            \mathcal{R}(\hat{\theta}_n) \lesssim \frac{ \|f^*\|_{\mathcal{B}}^2}{m} + \lambda \|f^*\|_{\mathcal{B}} \sqrt{\frac{\log (2d)}{n}} + \sqrt{\frac{\log(1/\delta) + \log(n)}{n}}.
        \]
\end{theorem}

Similar approximation theory and a priori error estimates have been proved for other machine learning models.
Here is a brief summary of these results.

\bi
\item Random feature model:  The corresponding function space is the reproducing kernel Hilbert space (RKHS).
\item Residual networks (ResNets):  The corresponding function space is the so-called flow-induced space
introduced in \cite{EMaWu2019}.
\item Multi-layer neural networks:  A candidate for the corresponding function space is the multi-layer space
 introduced in \cite{multilayerspace}.
\ei
What is really important is the ``norms'' that control the approximation error and the generalization gap.
These quantities are defined for functions in the corresponding spaces.
After the approximation theorems and Rademacher complexity estimates are in place, one can readily prove
a theorem of the following type for regularized models:
 Up to logarithmic terms,  the minimizers of the regularized models satisfy:
$$
\mathcal{R}(\hat{f}) \lesssim \frac{\Gamma(f^*)}{m} + \frac{\gamma(f^*)}{\sqrt{n}}
$$
where $m$ is the number of free parameters, $n$ is the  size of the training dataset.
Note that for the multilayer spaces, the results proved in \cite{multilayerspace} are not as sharp.

We only discussed the analysis of the hypothesis space. There are many more other questions. 
We refer to \cite{Review-3} for more discussion on the current understanding of neural network-based machine learning.

\section{Machine learning from a continuous viewpoint }

Now we turn to alternative formulations of machine learning. Motivated by the situation for PDEs, we would like to first formulate
machine learning in a continuous setting and  
then discretize to get concrete models and algorithms.  The key here is that continuous problems that we come up with should be
 nice mathematical problems. For PDEs, this is accomplished by requiring them to be ``well-posed''.
 For problems in calculus of variations,  we require the problem to be ``convex'' in some sense and lower semi-continuous.
 The point of these requirements is to make sure that the problem has a unique solution.
Intuitively, for machine learning problems, being  ``nice'' means that the  variational problem should  have a simple landscape.
How to formulate this precisely is an important research problem for the future.

As was pointed out in \cite{EMaWucontinuous}, the key ingredients for the continuous formulation are as follows:
\bi
\item representation of functions (as expectations)
\item formulating the variational problem (as expectations)
\item optimization, e.g. gradient flows
\ei

\subsection{Representation of functions}

Two kinds of representations are considered in \cite{EMaWucontinuous}: Integral transform-based representation and
flow-based representation.
The simplest integral-transform based representation is a generalization of \eqref{Fourier-alternative}:
\begin{align*}
f(\bx; \theta) = & \int_{\R^d} a(\wb) \sigma (\wb^T \bx) \pi(d\wb)\\
 =& \E_{\wb \sim \pi} a(\wb) \sigma (\wb^T\bx) \\
 = &\E_{(a, \wb) \sim \rho} a \sigma (\wb^T\bx) \\
 = & \E_{\bu \sim \rho} \phi (\bx, \bu)
\end{align*}
Here $\theta  $ denotes the parameters in the model: 
$\theta$ can be $a(\cdot)$ or the prob distributions $\pi$ or $\rho$.

This representation corresponds to two-layer neural networks.  A generalization to multi-layer neural networks is 
presented in \cite{multilayerspace}.

Next we turn to flow-based representation:
\begin{align}
 \frac{d \bz}{d\tau} = &\E_{\bw \sim \pi_\tau} \ba(\bw, \tau) \sigma(\bw^T \bz)\\
= &\E_{(\ba, \bw) \sim \rho_\tau}  \ba  \sigma(\bw^T \bz) \\
= & \E_{\bu \sim \rho_\tau} \phi (\bz , \bu), \quad \bz(0, \bx)= \bx 
\end{align}
\[
f(\bx, \theta) = \bm{1}^T \bz(1, \bx)
\]
In this representation, the parameter
$\theta$ can be either $ \{a_\tau(\cdot)\}$ or  $\{\pi_\tau\}$ or $\{\rho_\tau\}$

\subsection{The stochastic optimization problem}

Stochastic optimization problems are of the type:
  \[
  \min_\theta \, \EE_{\bw \sim \nu} g(\theta, \bw)
  \]
  These kinds of problems can readily  be approached by stochastic algorithms, which is a key component in machine learning.
  For example, instead of the gradient descent algorithm:
  $$\theta_{k+1} = \theta_k - \eta \nabla_{\theta} \EE_{\bw \sim \nu} g(\theta, \bw)
=\theta_k - \eta \nabla_{\theta} \sum_{j=1}^n g(\theta, \bw_j)
$$
one can use the stochastic gradient descent:
$$\theta_{k+1} = \theta_k - \eta \nabla_{\theta} g(\theta, \bw_k)
$$
where $\{\bw_k\}$ is a set of random variables sampled from $\nu$.

The following are some examples of the stochastic optimization problems that arise in modern machine learning:
\bi
\item Supervised learning:  In this case, the minimization of the population risk becomes
$$
\mathcal{R}(f) = \EE_{\bx \sim \mu} (f(\bx) - f^*(\bx))^2
 $$
\item Eigenvalue problems for quantum many-body Hamiltonian:
$$
 I(\phi)  =  \frac{(\phi, H\phi)}{(\phi, \phi)} = \EE_{\bx \sim \mu_{\phi}} \frac{\phi(\bx) H \phi(\bx)}{\phi(\bx)^2}, \quad  \mu_\phi (d\bx)= \frac 1Z |\phi(\bx)|^2  d\bx
$$
Here $H$ is the Hamiltonian of the quantum system.
\item Stochastic control problems:
$$
  L({\{a_t\}_{t =0}^{T-1}}) = \mathbb{E}\big\{\sum_{t =0}^{T-1} c_t(s_t, a_t))+c_T(s_T)\big\}
  $$
  \ei
  
  Substituting the representations discussed earlier to these expressions for the stochastic optimization problems, 
  we obtain the final variational problem that we need to solve.
  
  One can either discretize these variational problems directly and then solve the discretized problem using some optimization algorithms,
  or one can write down continuous forms of some optimization algorithms, typically gradient flow dynamics, and then discretize these
  continuous flows.  We are going to discuss the second approach.
  
 \subsection{Optimization: Gradient flows}

To write continuous form of the gradient flows, we draw some inspiration from statistical physics.  Take the supervised learning as 
an example. We regard the population risk as being the 
``free energy'', and following Halperin and Hohenberg \cite{HH1977}, we divide the parameters into two kinds,  conserved and non-conserved. For example,
 $a$ is a non-conserved parameter and $\pi$ is conserved since its total integral has to be 1.
 
 For non-conserved parameter, as was suggested in \cite{HH1977}, one can use the  ``model A'' dynamics:
\[
\frac{\partial a}{\partial t} = - \frac{\delta \CR}{\delta a}
\]
which is simply the usual $L^2$ gradient flow.

For conserved parameters such as $\pi$, one should use the ``model B'' dynamics which works as follows:
First define the ``chemical potential'' 
\[
V= \frac{\delta \CR}{\delta \pi}.
\] 
From the chemical potential, one obtains the velocity field $\vb$ and the current $J$:
\[
{\bf J}= \pi \vb, \, \vb= - \nabla V
\]
The continuity equation then gives us the gradient flow dynamics:
\[
\frac{\partial \pi}{\partial t} + \nabla \cdot {\bf J} = 0.
\]
This is also the gradient flow under the Wasserstein metric.

\subsection{Discretizing the gradient flows}

To obtain practical models, one needs to discretize these continuous problems. The first step is to
replace the population risk by the empirical risk using the training data.

The more non-trivial issue is how to discretize the gradient flows in the parameter space.
The parameter space has the following characteristics: (1) It has a simple geometry -- unlike the real space
which may have a complicated geometry.  (2) It is also usually high dimensional.
For these reasons, the most natural numerical methods for the  discretization in the parameter space is the 
particle method which is the dynamic version of Monte Carlo. Smoothed particle method might be helpful to improve
the performance. In relatively low dimensions, one might also consider the spectral method, particularly some sparse
version of the spectral method, due to the relative simple geometry of the parameter space.


Take for example the discretization of the conservative flow for the integral transform-based representation.
With the representation: $ f(\bx; \theta) = \E_{(a, \wb) \sim \rho} a \sigma (\wb^T\bx) $,  the gradient flow equation becomes:
\beq
\partial_{t} \rho = \nabla (\rho \nabla V),  \quad V= \frac{\delta \mathcal{R}}{\delta \rho}
\label{conserv-GD}
\eeq
The particle method discretization is based on:
$$\rho(a,  \wb, t) \sim \frac 1 m \sum_j \delta_{(a_j(t), \wb_j(t))} = \frac 1 m \sum_j \delta_{\bu_j(t)}
$$ 
where $\bu_j(t) = (a_j(t), \wb(t))$.
One can show that in this case, \eqref{conserv-GD} reduces to
$$
\frac{d \bu_j}{dt} = - \nabla_{\bu_j} I(\bu_1, \cdots, \bu_m)
$$
where
$$
I(\bu_1, \cdots, \bu_m) = \mathcal{R}(f_m), \quad
\bu_j = (a_j, \bw_j),  \quad f_m(\bx) = \frac 1m \sum_j  a_j \sigma(\wb_j^T \xb)
$$
This is exactly  gradient descent for ``scaled'' two-layer neural networks.
In this case, the continuous formulation also coincides with
the ``mean-field''  formulation for two-layer neural networks
\cite{ChizatBach,Mei,Eric,SS}.

\begin{figure}
\centering
\includegraphics[width=0.345\textwidth]{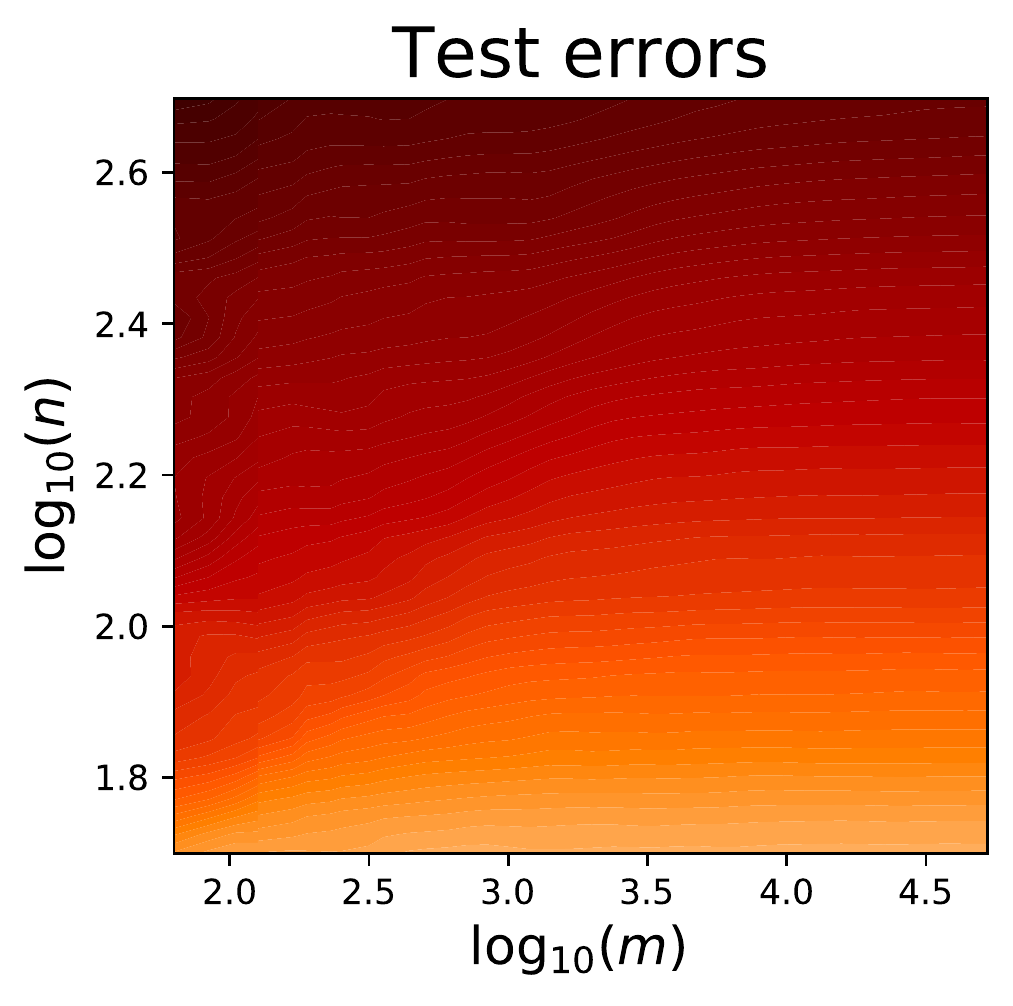}
\hspace*{0.5cm}
\includegraphics[width=0.4\textwidth]{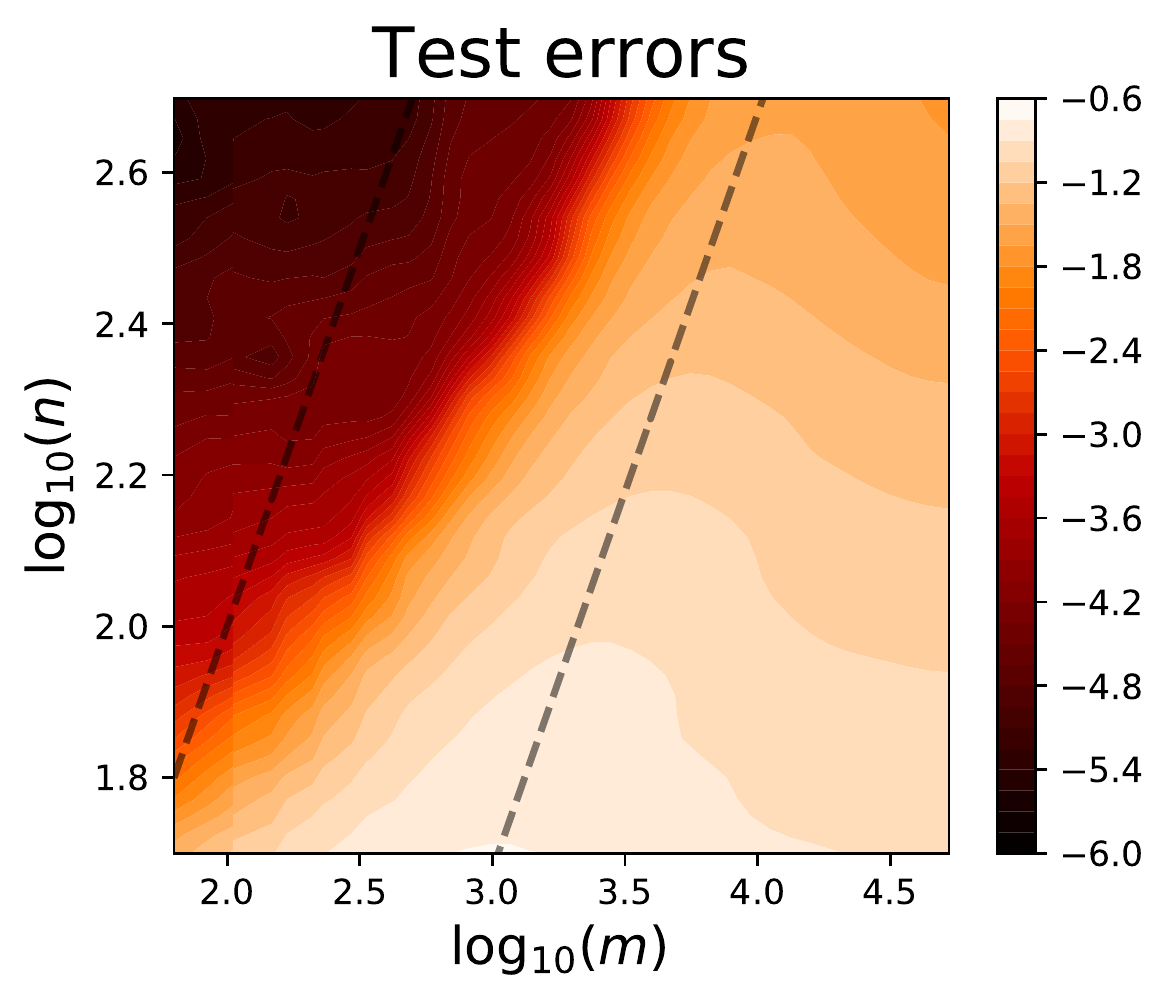}
\caption{\textbf{(Left)} continuous viewpoint; \textbf{(Right)} conventional NN models. Target function is a single neuron.
Reproduced with permission from Lei Wu. }
\label{heatmap}
\end{figure}

The scaling factor $1/m$ in front of the ``scaled'' two-layer neural networks is actually quite important and
makes a big difference for the test performance of the network.
Shown in Figure \ref{heatmap} is the heat map of the test error for two-layer neural network models with
and without this scaling factor.  The target function is the simple single neuron function: $f^*(\bx) = \sigma (x_1)$.
 The important observation is that in the absence of this scaling factor, the
test error shows a ``phase transition'' between a phase where the neural network model performs like an associated random 
feature model and another phase where it shows much better performance than the random feature model.
Such a phase transition is one of the reasons that choosing the right set of hyper-parameters, here the network width
$m$, is so important. 
However, if one uses the scaled form, this phase transition phenomenon is avoided and the performance is more
robust \cite{MaWuE2020}.

\subsection{The optimal control problem for flow-based representation}

The flow-based representation naturally leads to a control problem.
This viewpoint has been used explicitly or implicitly in machine learning for quite some time (see for example \cite{LeCun}).
The back propagation algorithm, for example, is an example of control-theory based algorithm.
Another more recent example is the development of maximum principle-based algorithm,  first introduced
in \cite{Qianxiao}.
Despite these successes, we feel that there is still a lot of room for using the control theory viewpoint to develop new algorithms.

We consider the flow-based representation in a slightly more general form
\[
 \frac{d \bz}{d\tau}
= \E_{\bu \sim \rho_\tau}  \bm{ \phi} (\bz, \bu), \quad \bz(0, \bx) = \bx
\]
where $\bz$ is the state, $ \rho_\tau$ is the control at time $\tau$.
Our  objective  is to   minimize $\mathcal{R}$ over $\{\rho_\tau\}$
\[
\mathcal{R}(\{\rho_\tau\}) = \EE_{\bx \sim \mu} (f(\bx) - f^*(\bx))^2  = 
\int_{\RR^d} (f(\bx) - f^*(\bx))^2 d \mu
\]
where as before 
\begin{equation}
f(\bx) = \bm{1}^T \bz(1, {\bx})
\end{equation}

One most important result for this control problem is Pontryagin's maximum principle (PMP).  To state this result, let us
define the Hamiltonian $H:\RR^d\times \RR^d \times \mathcal{P}_2(\Omega): \mapsto \RR$ as
\[
H(\bz,\bm{p},\mu) = \EE_{\bu\sim\mu}[\bm{p}^T\bm{\phi}(\bz,\bu)].
\]
Pontryagin's maximum principle asserts that the solutions of the control problem must satisfy:
\begin{align}
\rho_\tau = \mbox{argmax}_{\rho} \, \EE_{\bx}[ H\left(\bz^{t,\bx}_\tau, \bm{p}^{t,\bx}_\tau,\rho\right)],
 \,\, \forall \tau \in [0,1],
 \label{PMP}
\end{align}
where for each $\bx$, $\{(\bz_\tau^{t,\bx}, \bm{p}_\tau^{t,\bx})\}$ are defined by the forward/backward equations:
\begin{equation}\label{eqn: forward-back}
\begin{aligned}
\frac{d\bz_\tau^{t,\bx}}{d\tau} &= \nabla_{\bm{p}} H =  \EE_{\bu\sim\rho_\tau(\cdot;t)}[\bm{\phi}(\bz_\tau^{t,\bx},\bu)] \\
\frac{d\bm{p}_\tau^{t,\bx}}{d\tau} &= - \nabla_{\bz} H = \EE_{\bu\sim\rho_\tau(\cdot;t)}[\nabla_{\bz}^T\bm{\phi}(\bz_\tau^{t,\bx},\bu)\bm{p}^{t,\bx}_\tau],
\end{aligned}
\end{equation}
with the boundary conditions:
\begin{align}
    \bz_0^{t,\bx} &= {\bx} \\
    \bm{p}_1^{t,\bx} &= - 2(f(\bx;\rho(\cdot;t))-f^*(\bx)) \bm{1}. 
\end{align}

Pontryagin's maximum principle is slightly stronger than the KKT condition for the stationary points in that  \eqref{PMP}  is a statement
of optimality rather than criticality. In fact \eqref{PMP} also holds when the parameters are discrete and this has been used in \cite{Qianxiao2}
to develop efficient numerical algorithms for this case.

With the help of PMP, it is also easy to write down the gradient descent flow for the optimization problem.
Formally, one can simply write down the gradient descent flow for \eqref{PMP} for each $\tau$:
\begin{align}\label{eqn: flow-resnet}
\partial_{t} \rho_\tau(\bu,t) = \nabla \cdot \left(\rho_\tau(\bu,t)\nabla V(\bu;\rho)\right), \,\, \forall \tau \in [0,1],
\end{align}
where 
\[
V(\bu;\rho) = \EE_{\bx}[\frac{\delta H}{\delta \rho}\left(\bz^{t,\bx}_\tau, \bm{p}^{t,\bx}_\tau,\rho_\tau(\cdot;t)\right)],
\]
and $\{(\bz_\tau^{t,\bx}, \bm{p}_\tau^{t,\bx})\}$ are defined as before by the forward/backward equations.

To discretize the gradient flow, we can simply use:
\begin{itemize}
\item forward Euler for the flow in $\tau$ variable, with step size $1/L$ ($L$ is the number of grid points in the $\tau$ variable);
\item particle method for the gradient descent dynamics, with $M$ samples in each $\tau$-grid point.
\end{itemize}
This gives us
\begin{align}
\bz_{l+1}^{t,\bx} &= \bz_l^{t,\bx} + \frac{1}{LM}\sum_{j=1}^M \bm{\phi}(\bz^{t,\bx}_l,\bu^{j}_{l}(t)) ,\quad l=0,\dots, L-1\\
\bm{p}_{l}^{t,\bx} &= \bm{p}_{l+1}^{t,\bx} + \frac{1}{LM}\sum_{j=1}^M \nabla_{\bz} \bm{\phi}(\bz_{l+1}^{t,\bx},\bu_{l+1}^j(t)) \bm{p}_{l+1}^{t,\bx} ,\quad l =0,\dots, L-1\\
\frac{d \bu_l^{j}(t)}{d t} &=  -\EE_{\bx}[\nabla_{\bw}^T\bm{\phi}(\bz_l^{t,\bx},\bu_l^j(t)) \bm{p}_l^{t,\bx}].
\end{align}
This recovers the GD algorithm (with back-propagation) for the (scaled) ResNet:
\[
\bz_{l+1} = \bz_{l} + \frac{1}{LM}\sum_{j=1}^M \bm{\phi}(\bz_l,\bu_l).
\]
We call this ``scaled'' ResNet because of the presence of the factor $1/(LM)$.

In a similar spirit, one can also obtain an algorithm using PMP \cite{Qianxiao}.
Adopting the terminology in control theory,   this kind of algorithms  are called ``method of successive approximation'' (MSA).
The basic MSA is as follows:
\bi     
\item      Initialize: $\theta^0\in \mathcal{U}$ 
\item For $k = 0, 1, 2, \cdots$:
        \bi
        \item  Solve 
        $$\frac{d {\bz}^{k}_\tau}{d \tau} = \nabla_{\bm{p}} H(\bz^{k}_\tau,{\bm{p}}^{k}_\tau,\theta^k_\tau), \quad \bz^{k}_0 =  \bx $$
        \item  Solve 
        $$ \frac{d \bm{p}^{k}_\tau}{d \tau}= -\nabla_{\bz} H(\bz^{k}_\tau,{\bm{p}}^{k}_\tau,\theta^k_\tau), \quad \bm{p}^{k}_1 = - 2(f(\bx;\theta^k)-f^*(\bx)) \bm{1}$$
        \item  Set 
        $$\theta^{k+1}_\tau = \argmax_{\theta\in\Theta} H(\bz^{k}_\tau,\bm{p}^{k}_\tau,\theta)$$
         for each $ \tau \in[0,1]$
        \ei
\ei
In practice, this basic version does not perform as well as the  ``extended MSA''  which works in the same way as the MSA except that the Hamiltonian is replaced by an extended Hamiltonian \cite{Qianxiao}:
\[
\tilde{H}(\bz,\bm{p},\theta, \bv, \bm{q}) :=H(\bz,\bm{p},\theta) - \frac{1}{2} \lambda \Vert \bv-f(\bz, \theta) \Vert^2
 - \frac{1}{2} \lambda \Vert \bm{q}+\nabla_{\bz} H(\bz, \bm{p},\theta) \Vert^2.
\]

\begin{figure}[H]
	\centering
	\includegraphics[width=0.9\textwidth]{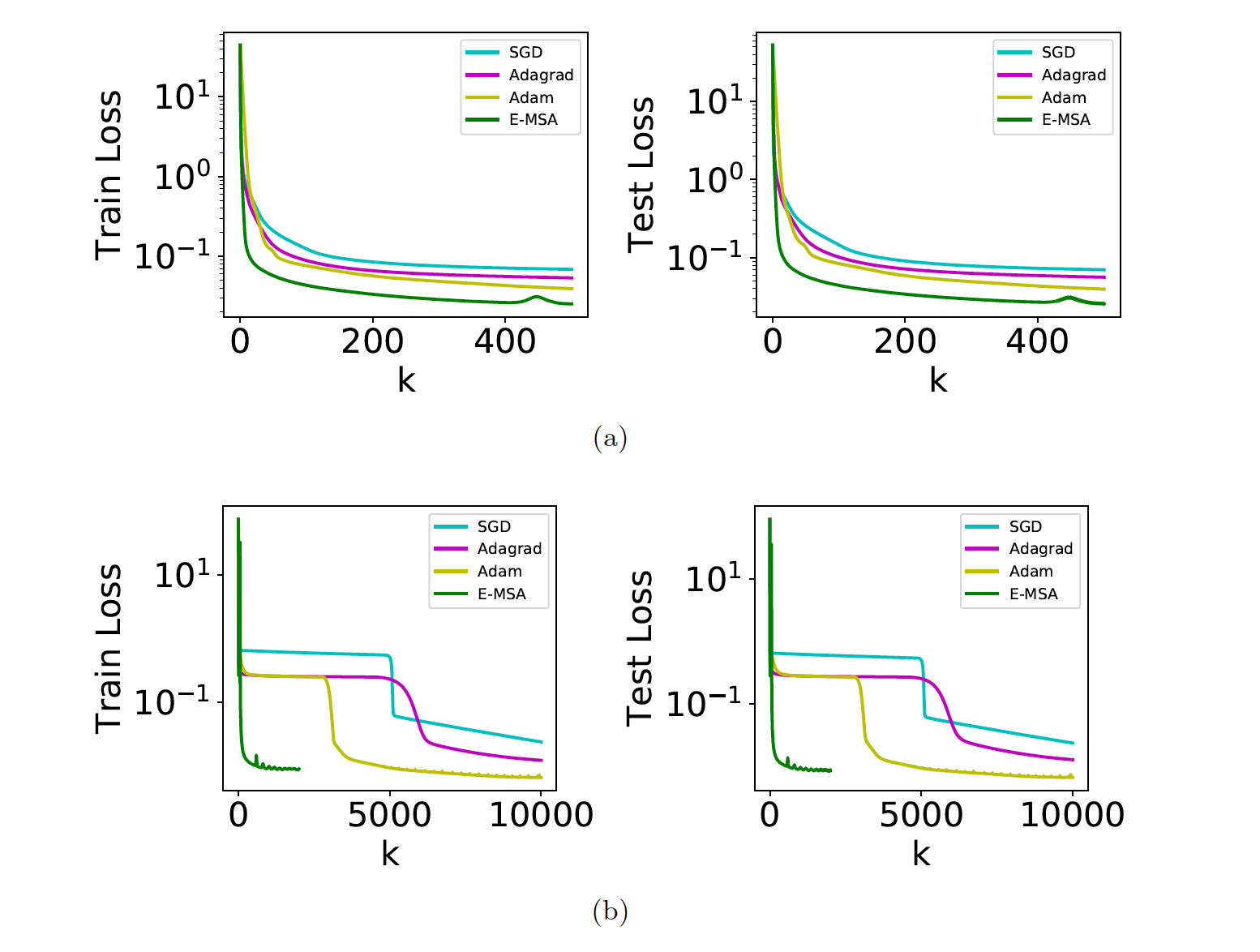}
	\vspace{-1em}
	\caption{Comparison of the extended MSA with different versions of stochastic gradient descent algorithms.
	The top figures show results for small initialization.  The bottom figures show results for bigger initialization.
	Reproduced with permission from Qianxiao Li.  See also \cite{Qianxiao}. }
  \label{MSA}
\end{figure}

Figure \ref{MSA} shows the results of the extended MSA compared with different versions of SGD for two
kinds of initialization. One can see that in terms of the number of iterations, extended MSA outperforms all the SGDs.
In terms of wall clock time, the advantage of the extended MSA is diminished significantly.  This is possibly
due to the inefficiencies in the implementation of the optimization algorithm (here the BFGS) used for solving  \eqref{PMP}.
We refer to \cite{Qianxiao} for more details.  In any case, it is clear that there is a lot of room for improvement.

\section{Concluding remarks}

In conclusion, we have discussed a wide range of problems for which machine learning-based algorithms have made and/or
will make a significant difference. These problems are relatively new to computational mathematics.
We believe strongly that machine learning-based algorithms will also significantly impact the way we solve more traditional
problems in computational mathematics. However,  research in this direction is still at a very early stage.  

Another important area that machine learning can be of great help is multi-scale modeling. The moment-closure problem discussed above
is an example in this direction. There are many more possible applications, see \cite{E2011}.
Machine learning seems to be able to provide the missing link in making advanced  multi-scale modeling techniques really
practical.  For example in the heterogeneous multi-scale method (HMM) \cite{Acta, EE03}, one important component is to 
extract the relevant macro-scale information from micro-scale simulation data. This step has always been a major obstacle in HMM.
It seems quite clear that machine learning techniques can of great help here.

We also discussed how the viewpoint of numerical analysis can help to improve
the mathematical foundation of machine learning as well as propose new and possibly more robust formulations.
In particular, we have given a taste of how high dimensional approximation theory should look like.
We also demonstrated that commonly used machine learning models and training algorithms can be recovered from
some particular discretization of continuous models, in a scaled form.
From this discussion, one can see that neural network models are quite natural and rather inevitable. 

What have we really learned from machine learning?
Well,  it seems that the most important new insight from machine learning is the representation of functions as expectations.
We reproduce them here for convenience:
\bi
\item integral-transform based:
\[
 f(\bx) =  \E_{(a, \bw) \sim \rho} a \sigma (\bw^T\bx)
\]
\[
         f(\bx)= \E_{\theta_L \sim \pi_L}  a^{(L)}_{\theta_L} \sigma(\E_{\theta_{L-1} \sim \pi_{L-1}}\dots
     \sigma(\E_{\theta_1 \sim \pi_1} a^1_{\theta_2,\theta_1}\sigma( 
         a^0_{\theta_1}\cdot  \bx) )\dots  )
\]
\item flow-based:
\begin{align}
 \frac{d \bz}{d\tau} = &\E_{(\ba, \bw) \sim \rho_\tau}  \ba  \sigma(\bw^T \bz), \quad \bz(0, \bx)= \bx  \\
f(\bx, \theta) = & \bm{1}^T \bz(1, \bx)
\end{align}
\ei

From the viewpoint of computational mathematics, this suggests that the central issue will move from specific discretization
schemes to more effective representations of functions.

This review is rather sketchy. Interested reader can consult the  three review articles \cite{Review-1, Review-2, Review-3}
for more details.

{\bf   Acknowledgement:} I am very grateful to my collaborators for their contribution to the work described here.
In particular, I would like to express my sincere gratitude to Roberto Car,  Jiequn Han,  Arnulf Jentzen, Qianxiao Li, Chao Ma, Han Wang, 
Stephan Wojtowytsch, and Lei Wu for the many discussions that we have had on the issues discussed here.
This work is supported in part by a gift to the Princeton University from iFlytek as well as the ONR grant N00014-13-1-0338.

\end{document}